\documentclass[final]{siamltex}
\usepackage[T1]{fontenc}
\usepackage{graphicx}
\usepackage{amssymb,latexsym,amsmath}
\usepackage{mathrsfs}
\usepackage{enumitem}
\usepackage{caption}
\usepackage{algorithm}
\usepackage{algpseudocode}
\usepackage{comment}
\usepackage{inputenc}
 \usepackage{makeidx}
\usepackage{color}
\usepackage{graphicx}
\graphicspath{ {./imgs/} }
\usepackage{babel}
\usepackage{cite}
\usepackage{amssymb,latexsym,amsmath,mathtools}
\usepackage{xcolor}
\usepackage{textcomp}
\usepackage{verbatim}
\usepackage{float}
\usepackage{titlesec}
\usepackage[mathlines]{lineno}

\usepackage{tikz}
\usetikzlibrary{shapes, arrows, positioning}

\titleformat{\paragraph}
{\normalfont\normalsize\bfseries}{\theparagraph}{1em}{}
\titlespacing*{\paragraph}
{0pt}{0.25ex plus 1ex minus .2ex}{1.5ex plus .2ex}

\allowdisplaybreaks

\newtheorem{assumption}[theorem]{Assumption}

\newtheorem{remark}[theorem]{Remark}
\DeclareMathOperator*{\argmin}{argmin}

\def\qed{\hfill $\diamond$}

\allowdisplaybreaks

\renewcommand{\theenumi}{\roman{enumi}}

\title{\bf Kernel Metrics and Learning for Borel MDPs: Identifiability and Adaptive Control}
\author{Omar Mrani-Zentar$^{1}$,  and Serdar Y\"uksel$^{1}$
\thanks{$^{1}$Are with the Department of Mathematics and Statistics at Queen's University, Kingston ON, Canada.
{\tt\small o.mranizentar@queensu.ca corresp. yuksel@queensu.ca}}%
}
\begin{document}

\maketitle
\begin{abstract}
     We consider a Markov decision process with standard Borel spaces and an unknown transition kernel under the average cost criterion. We do not impose any parametrization on the set of possible kernels. To facilitate our analysis, we first develop implication relations between several topologies on kernels defined by pointwise, continuous, or uniform weak convergence; we then review robustness properties on the space of kernels, and finally we establish compactness conditions on the space of kernels. Building on this regularity analysis, we then present two data-driven identifiability results; the first one being Bayesian and the second one empirical. Our conditions for the Bayesian setting are significantly more relaxed compared with prior work which considered either finite or parametric models, though we do not obtain a rate of convergence. Our analysis is asymptotic and builds on measurability in terms of the tail $\sigma$-field of the available information. Identifiability results are then used to design near-optimal adaptive control policies which alternate between periods of exploration, where the controller acts according to a policy which is conducive to the identification of the true kernel, and periods of exploitation where the controller's information on the kernel is utilized. We will establish that such policies are near optimal. In summary, our contribution is with regard to the general standard Borel setup where there is no apriori parametric representation.  
\end{abstract}
\section{Introduction}
In many practical stochastic control settings, one does not know all the parameters involved in the problem. For instance, the controller may not know the functional form of the dynamics or the stochastic kernel driving the state process \cite{Hernandez-lerma1987AdaptiveControl,Adelman2024Thompsonsampling,Kim2017TSfinite-case-average-cost,Kim2019TS-continuous-case-average-cost,Caines-adaptive-control, On-the-certainty-equivalence,Sznaier2025DataDrivenLearningControl,CostaDufourGenadot2026}. In other settings, the controller may know the function through which the state evolution is realized but may not know the distribution of the driving i.i.d noise \cite{Hernandez-Lerma1990DensityestimationAndAdaptiveControl,Gordienko1985}. For such setups where the model is not known, learning theoretic methods have been studied. These involve (i) Bayesian (see e.g. \cite{Kim2017TSfinite-case-average-cost,Kim2019TS-continuous-case-average-cost}), (ii) empirical (see e.g. \cite{Gordienko1985,DuPr14,Hernandez-Lerma1990DensityestimationAndAdaptiveControl,yichen2024}), and (iii) reinforcement learning methods (see e.g. \cite{gosavi2004reinforcement,abounadi2001learning,ky2026qaverage}). Each of these methods have their respective attributes; notably, Bayesian learning is more effective when there is limited data, empirical learning requires more lenient regularity conditions, and reinforcement learning does not require access to sophisicated mathematical tools and knowledge for implementation. 

It is important to note that empirical learning and reinforcement learning methods have an intrinsic equivalence connection, see e.g. \cite[Section 5.3]{ky2026qaverage}. Regarding empirical learning, by utilizing the convergence of the empirical processes in the Wasserstein sense, \cite{Gordienko1985,DuPr14,yichen2024} have established consistency results. For general weakly continuous models, this has been studied in \cite{kara2022robustness} under almost sure weak convergence of empirical estimates. If the driving noise is known to be absolutely continuous with respect to the Lebesgue measure, then empirical measures converge in total variation, and stronger results are obtained \cite{Hernandez-Lerma1990DensityestimationAndAdaptiveControl} \cite{kara2022robustness}.


In contrast to empirical or reinforcement learning, Bayesian learning methods have been applied to models which are either finite or continuous models which admit a parametric representation. One popular approach for Bayesian learning is via Thompson sampling \cite{russo2018tutorial,Adelman2024Thompsonsampling}. This approach, in the context of model learning, typically requires absolute continuity of all the kernels with respect to the same reference measure as well as a parametrization of the Radon-Nikodym derivates. Thompson sampling consists of sampling the parameter from a prior distribution, and then applying an optimal policy that corresponds to the sampled parameter. Subsequently, the new state and reward accrued are observed and that information is used to perform a Bayesian update on the prior. For finite and parametric continuous space MDPs, respectively, \cite{Kim2017TSfinite-case-average-cost,Kim2019TS-continuous-case-average-cost} developed asymptotic convergence results in the regret framework by imposing passive learning, similar to \cite{Adelman2024Thompsonsampling}, and Lipschitz continuity of densities. In \cite{bayrooti2025noregret}, an alternative Bayesian adaptive control algorithm is presented and convergence rates are established when the transition kernel and rewards are modeled as Gaussian.

In this paper, we will consider an average cost problem where the state and action spaces are standard Borel and the transition kernel is unknown. A parametric representation is not assumed. The generality of the setup requires that one develops kernel theoretic mathematical analysis tools. 

In the context of such a kernel based analysis, under the average cost criterion, the continuity of the optimal cost in the kernel was established under the uniform total variation metric in \cite{Hernandez-lerma1987AdaptiveControl}. Additionally, \cite{Hernandez-lerma1987AdaptiveControl} shows that for any strongly consistent sequence of parameters, i.e., a sequence of estimators that converge almost surely to the true parameter, one can construct an adaptive policy which is optimal. Analogous results on continuity and robustness have been obtained in \cite{kara2022robustness} (building on \cite{Lan81,muller1997does} which considered continuity of value functions and \cite{kara2020robustness} which addressed robustness and stability) under continuous weak convergence on kernels, which is a more practical condition for empirical learning based methods. To complete our literature review on kernel convergence, which is an emerging research topic, we note that a closely related topology on spaces of probability measures corresponding to laws of stochastic processes is defined by the following criterion: a sequence of stochastic processes converges to another process if their finite-dimensional marginals converge weakly, and their conditional distributions on future variables given the past (viewed as measure-valued stochastic processes) also converge weakly. \cite{aldous1981weak} has termed this {\it extended weak convergence} and \cite{hellwig1996sequential} {\it the information topology}; these have recently been shown to be equivalent in discrete-time \cite{backhoff2019all,pammer2024note}. The closely related {\it adapted Wasserstein metric} \cite{bartl2024wasserstein,backhoff2019all,beiglbock2022approximation} has been shown to possess strong robustness properties in a variety of applications \cite{bayraktar2020continuity,julio2020adapted,bartl2023sensitivity}, analogous to the role of continuous weak convergence discussed above. Further relevant work on a similar robustness framework includes \cite{GordienkoSystemControlLetters, GordienkoODE, GordienkoKybernetic,bozkurt2024modelapproximationmdpsunbounded} (discrete-time) and in \cite{pradhan2022robustness} (continuous-time). We refer the reader to \cite{saldi2025kernel} for further equivalence properties and an extensive literature review on various (strong and weak) topologies on stochastic kernels. 

In our paper, noting the critical dependence of robustness properties on the space of stochastic kernels and the particular notions of convergence under which learning is achieved, we will present several key supporting results establishing implication relations between various topologies on kernels. We also establish necessary compactness conditions on the space of kernels.

{\subsection{Notation}
   (i) For any standard Borel set $\mathbb{X}$ (Borel measurable subset of a complete and separable metric space) we will use the notation $\mathcal{B}(\mathbb{X})$ to denote the Borel sigma-algebra on $\mathbb{X}$. (ii) For any standard Borel set $\mathbb{X}$ we denote the space of probability measures on $\mathbb{X}$ under the weak convergence topology by $\mathcal{P}(\mathbb{X})$. (iii) We will denote a finite sequence $\{x_{t}\}_{t=l}^{t=k}$, where $k,l\in \mathbb{Z}^{+}$ by $x_{[l,k]}$. (iv) We will use $C_{b}(\mathbb{X})$ to denote the set of continuous and bounded functions on the space $\mathbb{X}$. (V) For a metric space $\mathbb{X}$ we will use the notation $d_{\mathbb{X}}(.,.)$ to denote its metric.}
\subsection{Problem Definition}


Consider the following general model.
\begin{eqnarray}\label{equationSystemchp2}
X_{t+1}=f(X_t,U_t,W_t),
\end{eqnarray}
where $X_t$ is an $\mathbb{X}$-valued state variable, $U_t$ a $\mathbb{U}$-valued control action variable, $W_t$ a $\mathbb{W}$-valued i.i.d noise process, and $f$ a measurable function. We assume that $\mathbb{X}, \mathbb{U},\mathbb{W}$ are {\it standard Borel}, i.e., Borel subsets of {\it Polish metric} spaces (complete and separable metric spaces). We assume that all random variables live in the probability space $(\Omega, {\cal F}, P)$. Using stochastic realization results (see \cite[Lemma 1.2]{gihman2012controlled}, or \cite[Lemma 3.1]{BorkarRealization}), it follows that the model above in (\ref{equationSystemchp2}) consists of all $(\mathbb{X} \times \mathbb{U})^{\mathbb{Z}_+}$-valued stochastic processes which satisfy the following characterization: for all Borel sets $B \in {\cal B}(\mathbb{X})$, $t \geq 0$, and $P$-almost all realizations $X_{[0,t]}, U_{[0,t]}$:
\begin{eqnarray} \label{eq_evol1}
P( X_{t+1} \in B | X_{[0,t]}=x_{[0,t]}, U_{[0,t]}=u_{[0,t]}) = P( X_{t+1} \in B | X_t=x_t, U_t=u_t) =: {\cal T}(B | x_t, u_t) \nonumber \\
\end{eqnarray}
where ${\cal T}(\cdot|x,u)$ is a {\it stochastic kernel} from $\mathbb{X} \times \mathbb{U}$ to $\mathbb{X}$ (so that for every $B$, ${\cal T}(B | \cdot, \cdot)
$ is a measurable function on $\mathbb{X} \times \mathbb{U}$, and for every fixed $(a, b) \in \mathbb{X} \times \mathbb{U}$, ${\cal T}(\cdot | a, b)$ is a probability measure on $(\mathbb{X}, {\cal B}(\mathbb{X})$). A stochastic process which satisfies (\ref{eq_evol1}) is called a {\it controlled Markov chain}. For a standard Borel space $\mathbb{X}$, we will use ${\cal P}(\mathbb{X})$ to denote the set of probability measures on $\mathbb{X}$.



We consider a problem where the controller at any $t \in \mathbb{Z}_+$ has access to past state and action realizations and the current state: $I_{t}=\{X_{[0,t]},U_{[0,t-1]}\}$. We consider the problem where the controller does not know the kernel $\mathcal{T}$ but aims to minimize an average cost criterion over all admissible policies $\gamma \in \Gamma$ where $\Gamma=\{\gamma=\{\gamma_{t}\}_{t\in\mathbb{Z}^{+}}\text{ with } \gamma_{t}:I_{t}\rightarrow \mathbb{U}\mid\gamma_{t}$ is $\sigma(I_{t})$-measurable$\}$. For later use, we further denote by $\Gamma_{S}$ the set of stationary Markovian control policies, which map the state to actions via a fixed map at each time stage. The cost criterion is given by  
\begin{equation} \label{cost criteria}
    J({\cal T},x,\gamma)=\limsup_{T\rightarrow \infty} E^{\gamma}_{x}[\frac{1}{T}\sum_{t=0}^{T-1}c(X_{t},U_{t})]
\end{equation}
where $x$ is the initial state, and $c:\mathbb{X}\times\mathbb{U}\rightarrow [0,\infty)$ is the cost functional. We require that $\mathbb{U}$ is compact. For the case of Bayesian learning, we assume that the controller also has a prior $\Theta$ on the space of transition kernels (which will admit a metric space structure) so that $\Theta\in \mathcal{P}(\mathbb{M})$ where $\mathbb{M}$ is the space of candidate transition kernels.

Throughout the paper, we will also make the following assumption
\begin{assumption}\label{AssumptionWFAverage} 
 \begin{enumerate}
    \item[(a)] The cost function $c: \mathbb{X} \times \mathbb{U} \to \mathbb{R}_+$ is continuous and bounded; that is, $c \in C_{b}(\mathbb{X} \times \mathbb{U})$. 
     \item[(b)] The transition kernel $\mathcal{T}$ is weak-Feller; that is, the map
    \begin{equation*}
        (x,u)\mapsto\int_{\mathbb{X}}f(z)\mathcal{T}(dz|x,u)
    \end{equation*}
    is continuous for all $f \in C_{b}(\mathbb{X})$. 
 \end{enumerate}
\end{assumption}

\subsection{Contributions}
\begin{itemize} 
\item {\bf Topologies on Kernels: Implication Relations, Conditions for Robustness and Continuity, and Compactness of the Space of Kernels.} To facilitate our kernel based analysis, we first develop relations between pointwise weak convergence, continuous weak convergence and uniform weak convergence. Notably, we show that pointwise weak convergence of kernels implies continuous weak convergence under the existence of equicontinuous densities (Proposition \ref{prop pointwise implies cont under densities}) and compactness of the set of candidate kernels under the bounded Lipschitz metric (Theorem \ref{kernel topologies}), respectively. These lead to the continuity of the optimal cost under pointwise weak convergence in kernels which, in turn, implies robustness of the optimal cost to kernel identifiability errors. Theorem \ref{them exisistence of eps-partitions} establishes an Arzel\'a-Ascoli type compactness criterion on the space of kernels. These mathematical results are critically utilized in the analysis to follow.

\item {\bf Bayesian Learning, Identifiability, and Adaptive Control.} For a controlled Markov model with standard Borel spaces, in Theorem \ref{martingale convergent}, we show that if an exploration policy, which leads to an ergodic process with an invariant probability measure which positively charges every non-empty open set, is adopted, then the true kernel can be learned up to any arbitrary error tolerance under uniform weak convergence through Bayesian updates. Moreover, we establish convergence of the max-likelihood estimator. The latter hinges critically on Lemma \ref{M is tail measurable} which shows that a random variable which represents the true kernel is measurable with respect to the tail $\sigma$-field. By Theorem \ref{kernel topologies}, which relates uniform weak convergence with continuous weak convergence, and utilizing robustness under this convergence, in Theorem \ref{adavptive control theorem} and \ref{martingale convergent} we construct a near optimal policy which alternates between periods of exploring the environment to learn the true kernel and periods of exploitation which utilize the knowledge accrued. Key to our argument is a measurability argument in terms of the tail $\sigma$-field of the information, and robustness to simultaneous experimentation and near optimal control utilizing continuous dependence of the expected cost under the Young topology, and Theorem \ref{exhaustive policies are dense} which shows that slight randomization leads to a small perturbation of policies under this topology.  

\item {\bf Empirical Learning, Identifiability, and Adaptive Control.} To complete our analysis, building on prior work, we present an algorithm (Algorithm \ref{alg:alg1}) that first quantizes the state and action spaces using the procedure suggested in \cite{SaLiYuSpringer} and then estimates the kernel of the resulting MDP through empirical occupation measures. In Theorem \ref{computational learning theorem I}, we establish the asymptotic convergence of the empirical adaptive learning algorithm to an optimal solution.
\end{itemize}

As can be seen from Figure \ref{fig:tikz-diagram} common to both learning paradigms is a uniform mixing condition (Assumption \ref{Assumption for ACOE}). Bayesian identifiability requires further assumptions in order to guarantee the compactness of the space of candidate kernels under the uniform bounded Lipschitz metric and for exploration to be sufficiently informative. Empirical learning on the other hand requires a minorization condition which leads to a suitable approximation of the original optimization problem via a model with finite state and action spaces.

\begin{figure}[t]
\centering
\begin{tikzpicture}[
    scale=0.6,
    transform shape,
    box/.style={
        rectangle, 
        draw, 
        minimum width=2.5cm, 
        minimum height=1cm, 
        align=center,
        fill=blue!10,
        rounded corners=2pt
    },
    arrow/.style={
        ->, 
        thick, 
        >=stealth
    },
    line/.style={
        thick
    },
    node distance=1.5cm and 1cm
]

\node[box] (top) {Measurable selection \\(Assumption \ref{AssumptionWFAverage})};

\node[box, below=of top] (middle) {Uniform mixing\\ (Assumption \ref{Assumption for ACOE})};

\draw[line] (top) -- (middle);

\node[box, below left=of middle] (left1) {Compactness of \\the space of candidate \\kernels (Assumption \ref{compactness under unif BL})};
\node[box, below right=of middle] (right1) {Minorization assumption\\ (Assumption \ref{assumption for quantizing actions})};

\draw[line] (middle.south) -- ++(0,-0.3) -| (left1.north);
\draw[line] (middle.south) -- ++(0,-0.3) -| (right1.north);


\node[box, fill=blue!10, below=of left1] (left3) {Bayesian Identifiability};


\node[box, fill=blue!10, below=of left3] (left4) {Adaptive Bayesian\\ Learning and Control};

\draw[arrow] (left3) -- (left4);

\draw[arrow] (left1) -- (left3);

\node[box, fill=blue!10, below=of right1] (right2) {Empirical\\ Identifiability};

\node[box, fill=blue!10, below=of right2] (right3) {Adaptive Empirical\\ Learning and Control};

\draw[arrow] (right1) -- (right2);
\draw[arrow] (right2) -- (right3);

\end{tikzpicture}
\caption{Overview of the assumptions leading to adaptive Bayesian learning and control as well as adaptive empirical learning and control.}
\vspace{-1cm}
\label{fig:tikz-diagram}
\end{figure}

Thus, our contribution is to present general conditions ensuring identifiability and adaptive learning. While prior work has typically required either parametric models or restrictive assumptions in the continuous setup (see e.g.  \cite[Assumption 2]{Adelman2024Thompsonsampling}, \cite[Assumption 2, 6]{Kim2019TS-continuous-case-average-cost}), our generality comes at the expense of not having a rate of convergence or sample complexity analysis unlike \cite{Adelman2024Thompsonsampling,Kim2019TS-continuous-case-average-cost}. If one imposes stronger regularity, such as uniform Wasserstein continuity on the kernels and cost functions, the analysis here can be generalized. 

\section{Kernel Regularity Results: Relations between Kernel Convergence Notions, Compactness of Spaces of Kernels and Robustness}

In this section, we present several key regularity results on spaces of kernels which will be critical for our learning and adaptive control analysis to follow. We will first provide relations between pointwise, uniform, and continuous convergence of sequences of kernels. Notably, we will present sufficient conditions which ensure continuous convergence of kernels. We also establish conditions for compactness of a given space of kernels. This is followed by two key supporting results which are reviewed: on robustness under continuous convergence of kernels, and continuous dependence of invariant measures on control policy under Young topology.

\subsection{Kernels and convergence notions}

In the following, we define several convergence notions on the set of kernels $\mathbb{M}$. Let us recall that for a Polish metric space $\mathbb{X}$, a sequence $\left\{\mu_n, n \in \mathbb{N}\right\}$ $\subset$ $\mathcal{P}(\mathbb{X})$ is said to converge to $\mu \in \mathcal{P}(\mathbb{X})$ weakly if and only if $\int_{\mathbb{X}} f(x) \mu_n(d x) \rightarrow \int_{\mathbb{X}} f(x) \mu(d x)$ for every continuous and bounded $f: \mathbb{X} \rightarrow \mathbb{R}$. One important property of weak convergence is that the space of probability measures on a complete, separable, metric (Polish) space endowed with the topology of weak convergence is itself complete, separable, and metric. An example of such a metric is the bounded Lipschitz metric, which is defined for $\mu, \nu \in \mathcal{P}(\mathbb{X})$ as $\rho(\mu, \nu):=\sup _{\| f \|_{BL\leq 1}}\left|\int f d \mu-\int f d \nu\right|$
where
$
\|f\|_{B L}:=\|f\|_{\infty}+\sup _{x \neq y} \frac{|f(x)-f(y)|}{d(x, y)}
$ and $\displaystyle \|f\|_{\infty}=\sup _{x \in \mathbb{X}}|f(x)|$.

For probability measures $\mu, \nu \in \mathcal{P}(\mathbb{X})$, the total variation metric is given by
\begin{eqnarray*}
    \|\mu-\nu\|_{T V}&=&2 \sup _{B \in B(\mathbb{X})}|\mu(B)-\nu(B)|
    =\sup_{\|f\|_{\infty \leq 1}}\left|\int f(x) \mu(\mathrm{d} x)-\int f(x) \nu(\mathrm{d} x)\right|,
\end{eqnarray*}
where the supremum is taken over all measurable real functions $f:\mathbb{X}\rightarrow\mathbb{R}$ such that $\displaystyle \|f\|_{\infty}=\sup _{x \in \mathbb{X}}|f(x)| \leq$ 1. A sequence $\mu_{\mathrm{n}}$ is said to converge in total variation to $\mu \in \mathcal{P}(\mathbb{X})$ if and only if $\| \mu_n-\mu\|_{T V} \rightarrow 0$. Building on the above, we now consider convergence notions for stochastic transition kernels:
\begin{definition}[Pointwise weak convergence]
    We say that a sequence of kernels $\{\mathcal{T}_{n}\}_{n\in \mathbb{N}}\subseteq \mathbb{M}$ converges to some kernel $\mathcal{T}$ weakly and denote $\mathcal{T}_{n}\rightarrow \mathcal{T} \in \mathbb{M}$ weakly, if and only if, for all $x\in \mathbb{X}$ and $u\in \mathbb{U}$: $\mathcal{T}_{n}(.|x,u)\rightarrow \mathcal{T}(.|x,u)$ weakly (pointwise). 
\end{definition}

\begin{definition}[Continuous weak convergence]
    We say that a sequence of kernels $\{\mathcal{T}_{n}\}_{n\in \mathbb{N}}\subseteq \mathbb{M}$ converges to some kernel $\mathcal{T}$  weakly continuously, if and only if, for any sequence $(x_{n},x_{n})\rightarrow (x,u)$ we have that $\mathcal{T}_{n}(.|x_{n},u_{n})\rightarrow \mathcal{T}(.|x,u)$ weakly. 
\end{definition}

\begin{definition}[Uniform weak convergence]
    We say that a sequence of kernels $\{\mathcal{T}_{n}\}_{n\in \mathbb{N}}\subseteq \mathbb{M}$ converges to some kernel $\mathcal{T}$ under the uniform bounded Lipschitz metric if and only if $\displaystyle \rho^{\mathrm{unif}}_{BL}(\mathcal{T}_{n},\mathcal{T}):=\sup_{(x,u)\in \mathbb{X}\times\mathbb{U}}\rho(\mathcal{T}_{n}(.|x,u),\mathcal{T}(.|x,u))\rightarrow 0$ as $n\rightarrow \infty$. 
\end{definition}

{
\begin{definition}[Uniform weak convergence on compact subsets] \label{unif weak on compact}
    We say that a sequence of kernels $\{\mathcal{T}_{n}\}_{n\in \mathbb{N}}\subseteq \mathbb{M}$ converges to some kernel $\mathcal{T}$ under the topology of uniform weak convergence on compact subsets if and only if for any compact subset $C\subset \mathbb{X}\times \mathbb{U}$ $\displaystyle \sup_{(x,u)\in C}\rho(\mathcal{T}_{n}(.|x,u),\mathcal{T}(.|x,u))\rightarrow 0$ as $n\rightarrow \infty$. 
\end{definition}
}
The above have focused on weak convergence based metrics as these are relatively less demanding for many applications. Nonetheless, one can state analogous notions of convergence under other metrics such as total variation or Wasserstein metrics or convergence notions which are not characterized by a metric, such as setwise convergence. The following convergence notion will be used for the space of stationary control policies. As before, we take $\mathbb{U}$ to be compact.
\begin{definition}\label{definition young Top}[Young Topology on the space of stationary Markov policies $\Gamma_{S}$] 
    Consider the set of probability measures on $\mathbb{X}\times \mathbb{U}$ with fixed marginal $\psi$ on $\mathbb{X}$: \[S=\bigg\{ P\in \mathcal{P}(\mathbb{X}\times \mathbb{U})\mid P(dx,du)=P(du|x)\psi(dx)\bigg\}\] A sequence of policies $\{\gamma_{n}\}\subseteq \Gamma_{S}$ is said to converge under the Young topology at input $\psi$ (see, e.g. \cite{BorkarRealization,yuksel2023borkar,BalderYoungMeasureTheory}) to $\gamma \in \Gamma_{S}$ if and only if $\int g(x,u)\gamma_{n}(du|x)\psi(dx)\rightarrow \int g(x,u)\gamma(du|x)\psi(dx)$ for any continuous and bounded function $g: \mathbb{X}\times \mathbb{U}\rightarrow \mathbb{R}$. This topology leads to a convex and compact formulation (see \cite[Section 2]{BorkarRealization} \cite[Section 3]{Bor91}). 
\end{definition}
{
\begin{remark}
    It follows that, since the marginal on the measurement space is fixed, the convergence in Definition \ref{definition young Top} is equivalent to the requirement that $\int g(x,u)\gamma_{n}(du|x)\psi(dx)\rightarrow \int g(x,u)\gamma(du|x)\psi(dx)$ for any measurable function $g$ which is bounded and continuous in $u$ (\cite[Theorem 3.10]{Schal},  \cite[Theorem 2.5]{balder2001}); such functions are also called Carath\'eodory functions.
\end{remark}
}

\subsection{Relations between convergence notions: Sufficient conditions for continuous convergence of kernels and compactness}


We start with presenting several assumptions which will be used to develop implication relations.

\begin{assumption} \label{all is WF}
   For all $M\in\mathbb{M}$, $M$ is weak-Feller; that is, $M(.|x_{n},u_{n})\rightarrow M(.|x,u)$ weakly whenever $(x_n,u_n) \to (x,u)$.
\end{assumption}
{We note that since $\mathcal{T}\in \mathbb{M}$, it follows that Assumption \ref{all is WF} implies Assumption \ref{AssumptionWFAverage} (b)}.
\begin{assumption} \label{CompactStateActionSpaces}
$\mathbb{X}$ and $\mathbb{U}$ are compact.
\end{assumption}

\begin{assumption} \label{compactness under unif BL}
    The set of all possible kernels endowed with the uniform bounded Lipschitz metric: $(\mathbb{M},\rho^{\mathrm{unif}}_{BL})$ is compact.
\end{assumption} 
{Theorem \ref{them exisistence of eps-partitions} presents several conditions which imply that $\mathbb{M}$ is pre-compact and hence its closure satisfies Assumption \ref{compactness under unif BL}}.

\begin{theorem} \label{kernel topologies}
    Let $\{\mathcal{T}\}_{n}\subseteq \mathbb{M}$ and $\mathcal{T} \in \mathbb{M}$. 
    \begin{itemize}
           \item [(i)] Suppose Assumption \ref{all is WF} holds. Then, $\mathcal{T}_{n}\rightarrow \mathcal{T} $ weakly uniformly implies that $\mathcal{T}_{n}\rightarrow \mathcal{T} $ weakly continuously.
        \item[(ii)] Suppose Assumption \ref{CompactStateActionSpaces} and \ref{all is WF} hold. If, $\mathcal{T}_{n}\rightarrow \mathcal{T}$ weakly continuously, then $\mathcal{T}_{n}\rightarrow \mathcal{T}$ weakly uniformly.  
     
        \item  [(iii)]   Suppose Assumptions \ref{compactness under unif BL}, and \ref{all is WF} hold. Suppose $\mathcal{T}_{n}\rightarrow \mathcal{T}$ pointwise weakly. Then, $\mathcal{T}_{n}\rightarrow \mathcal{T}$ weakly uniformly. 
    \end{itemize}
\end{theorem}


\textbf{Proof.} [i] Suppose $\mathcal{T}_{n}\rightarrow\mathcal{T} $ weakly uniformly. Let $(x_{n},u_{n})\rightarrow (x,u)$. Then, 
\begin{align*}
    \rho(\mathcal{T}_{n}(.|x_{n},u_{n}),\mathcal{T}(.|x,u))&\leq \rho(\mathcal{T}_{n}(.|x_{n},u_{n}),\mathcal{T}(.|x_{n},u_{n}))+\rho(\mathcal{T}(.|x_{n},u_{n}),\mathcal{T}(.|x,u))\rightarrow 0,
\end{align*}
since the first term converges to zero by uniform weak convergence, and the second term by the weak Feller property of ${\cal T}$. \\

\noindent[ii]  Suppose $\mathcal{T}_{n}\not \rightarrow\mathcal{T} $ weakly uniformly. Then, there exists $\epsilon>0$ and a sequence $(x_{n},u_{n})$ such that for all $n$ $\rho(\mathcal{T}_{n}(.|x_{n},u_{n}),\mathcal{T}(.|x_{n},u_{n}))\geq \epsilon$. By compactness, there exists a subsequence such that $(x_{n_{k}},u_{n_{k}})\rightarrow (\Bar{x},\Bar{u})\in\mathbb{X}\times\mathbb{U}$ and for all $k \in \mathbb{N}$, by the triangle inequality, 
\begin{align*}
   & \rho(\mathcal{T}(.|\Bar{x},\Bar{u}),\mathcal{T}(.|x_{n_{k}},u_{n_{k}}))+\rho(\mathcal{T}_{n_{k}}(.|x_{n_{k}},u_{n_{k}}),\mathcal{T}(.|\Bar{x},\Bar{u})) \\
&  \quad \quad \quad \quad \geq \rho(\mathcal{T}_{n_{k}}(.|x_{n_{k}},u_{n_{k}}),\mathcal{T}(.|x_{n_{k}},u_{n_{k}}))\geq \epsilon
\end{align*}
Since the first term above can be made arbitrarily small, we get that 
\[\liminf_{k\rightarrow \infty} \rho(\mathcal{T}_{n_{k}}(.|x_{n_{k}},u_{n_{k}}),\mathcal{T}(.|\Bar{x},\Bar{u}))\geq \frac{\epsilon}{2} >0.\] Thus ${\cal T}_n$ does not converge weakly continuously. \\

\noindent[iii]  Consider any subsequence $\{\mathcal{T}_{n_{k}}\}$. By Assumption \ref{compactness under unif BL}, there exists a further subsequence $\{\mathcal{T}_{n_{k_{j}}}\}$ which converges to some $\mathcal{T}^{*}$ weakly uniformly. Since $\mathcal{T}_{n}\rightarrow \mathcal{T}$, it follows that $\mathcal{T}_{n_{k_{j}}}\rightarrow \mathcal{T}$, hence, $\mathcal{T}^{*}\equiv\mathcal{T}$. We have thus shown that any subsequence of $\{\mathcal{T}_{n}\}$ has a further subsequence which converges weakly uniformly to $\mathcal{T}$. Now suppose $\mathcal{T}_{n}\not \rightarrow \mathcal{T}$ weakly uniformly. Then, there exists a subsequence $\mathcal{T}_{n_{k}}$ and $\epsilon >0$ such that for all $k \in \mathbb{N}$ $\rho^{\mathrm{unif}}_{BL}(\mathcal{T},\mathcal{T}_{n_{k}})>\epsilon$. However, this is a contradiction since $\mathcal{T}_{n_{k}}$ can not have a subsequence which converges weakly uniformly to $\mathcal{T}$.  \qed
 
\begin{assumption} \label{assumption for proposition II} 
All kernels $M\in \mathbb{M}$ admit densities $\{f^{M}_{x,u}(y)\}_{M\in \mathbb{M}}$ with respect to some probability measure $\psi$ and all such densities are equicontinuous over $(x,u)$ uniformly over $y$. 
\end{assumption}

\begin{proposition} \label{prop pointwise implies cont under densities}    
  Suppose that Assumption \ref{assumption for proposition II} holds. If $\mathcal{T}_{n}\rightarrow \mathcal{T}$ weakly (pointwise), then $\mathcal{T}_{n}\rightarrow \mathcal{T}$ weakly continuously. 
\end{proposition}

\textbf{Proof.}
 Let $g\in\mathbf{C}_{b}(\mathbb{X})$ and consider a sequence $(x_{n},u_{n})\rightarrow (x,u)$.
\begin{align*}
 &\bigg|\int g(y)\mathcal{T}_{n}(dy|x_{n},u_{n})-\int g(y)\mathcal{T}(dy|x,u)\bigg| =\bigg|\int g(y) f^{\mathcal{T}_{n}}_{x_{n},u_{n}}(y)\psi(dy)-\int g(y) f_{x,u}^{\mathcal{T}}(y)\psi(dy)\bigg|\\
    &\leq \bigg|\int g(y) f^{\mathcal{T}_{n}}_{x_{n},u_{n}}(y)\psi(dy)-\int g(y) f^{\mathcal{T}_{n}}_{x,u}(y)\psi(dy)\bigg| +  \bigg|\int g(y) f^{\mathcal{T}_{n}}_{x,u}(y)\psi(dy)-\int g(y) f_{x,u}^{\mathcal{T}}(y)\psi(dy)\bigg|\\
& \leq \|g\|_{\infty}\int \sup_{M\in \mathbb{M}}\bigg|f^{M}_{x_{n},u_{n}}(y)-f^{M}_{x,u}(y)\bigg|\psi(dy) +  \bigg|\int g(y) f^{\mathcal{T}_{n}}_{x,u}(y)\psi(dy)-\int g(y) f_{x,u}^{\mathcal{T}}(y)\psi(dy)\bigg| 
\end{align*}
By Assumption \ref{assumption for proposition II}, it follows that $\sup_{M\in \mathbb{M}}\bigg|f^{M}_{x_{n},u_{n}}(y)-f^{M}_{x,u}(y)\bigg|\rightarrow 0$ (pointwise in $y$) and hence, we get that $\int \sup_{M\in \mathbb{M}}\bigg|f^{M}_{x_{n},u_{n}}(y)-f^{M}_{x,u}(y)\bigg|\psi(dy)\rightarrow 0$. Because $\mathcal{T}_{n}\rightarrow \mathcal{T}$ weakly (pointwise), we have that $ \bigg|\int g(y) f^{\mathcal{T}_{n}}_{x,u}(y)\psi(dy)-\int g(y) f_{x,u}^{\mathcal{T}}(y)\psi(dy)\bigg|\rightarrow 0$ by definition. 
 \qed

Next, we provide sufficient conditions for Assumption \ref{compactness under unif BL} to hold.

{
\begin{assumption} \label{equicontionuity}
    The family of kernels $\mathbb{M}$ is uniformly weakly equicontinuous that is for any $\epsilon>0$, there exists $\delta>0$ such that for all $(x,u),(x',u')\in \mathbb{X}\times\mathbb{U}$: $d((x,u),(x',u'))<\delta \implies \forall M \in \mathbb{M}$: $\rho(M(.|x,u),M(.|x',u'))\leq \epsilon$.
\end{assumption}
}
The following can be viewed as an Arzel\'a-Ascoli type result for spaces of kernels, leading to conditions for pre-compactness of a set of stochastic kernels. For a general review of Arzel\'a-Ascoli type theorems see \cite[Chapter 7]{munkres2000topology}. The proof of Theorem \ref{them exisistence of eps-partitions} presents a more tailored argument for the context of this paper.
\begin{theorem} \label{them exisistence of eps-partitions}
Suppose Assumption \ref{equicontionuity} and \ref{CompactStateActionSpaces} hold, then $(\mathbb{M},\rho^{\mathrm{unif}}_{BL})$ is pre-compact, i.e., its closure is compact. 
\end{theorem}

\textbf{Proof.} Let $\epsilon>0$. By assumption \ref{equicontionuity}, for all $(x,u)\in \mathbb{X}\times\mathbb{U}$ there exists $\delta(x,u)>0$ such that for all $(x',u')\in \mathbb{X}\times\mathbb{U}$: $d((x,u),(x',u'))<\delta(x,u) \implies \forall M \in \mathbb{M}$: $\rho(M(.|x,u),M(.|x',u'))\leq \frac{\epsilon}{4}$. The open balls $\{B_{\delta(x,u)}(x,u)\}$ form an open cover of $\mathbb{X}\times \mathbb{U}$. Thus, by the compactness of $\mathbb{X}$, $\mathbb{U}$, there exists a finite cover for $\mathbb{X}\times\mathbb{U}$ using finitely many balls $\{B_{\delta(\Bar{x}_{i},\Bar{u}_{i})}(\Bar{x}_{i},\Bar{u}_{i})|i\in\{1,...,K\}\}$. Consider  $\mathbb{M}(x,u)=\{M(.|x,u)\mid M\in \mathbb{M}\}\subseteq \mathcal{P}(\mathbb{X})$. Because $\mathbb{X}$ is compact, $\mathcal{P}(\mathbb{X})$ is weakly compact, it follows that for all $(x,u)$ $\mathbb{M}(x,u)$ is precompact. Hence $ \bigcup_{i\in\{1,...,K\}} \mathbb{M}(\Bar{x}_{i},\Bar{u}_{i})$ is precompact. It follows that there exists an open cover of $\bigcup_{i\in\{1,...,K\}} \mathbb{M}(\Bar{x}_{i},\Bar{u}_{i})$ using open balls $\{B_{\frac{\epsilon}{4}}(a_{j})\}_{j=1}^{N}$ where for all $j$, $a_{j}\in \mathcal{P}(\mathbb{X})$. Hence, for all $M\in\mathbb{M}$ there exists a mapping $\sigma_{M}: \{1,...,K\}\rightarrow \{1,...,N\}$ such that $\rho(M(.|\Bar{x}_{i},\Bar{u}_{i}),a_{\sigma_{M}(i)})<\frac{\epsilon}{4}$. Let $\mathbb{M}_{\sigma}=\{M\in\mathbb{M}|\sigma_{M}=\sigma\}$. Then, $\{\mathbb{M}_{\sigma}\}$ forms a finite partition of $\mathbb{M}$. Consider any set of representative elements from this partition $\{\tau_{\sigma}\}_{\sigma}$ such that $\tau_{\sigma}\in \mathbb{M}_{\sigma}$. Let $(x,u)\in \mathbb{X}\times\mathbb{U}$. Let $i\in\{1,...,K\}$ be such that $(x,u)\in B_{\delta(\Bar{x}_{i},\Bar{u}_{i})}(\Bar{x}_{i},\Bar{u}_{i})$  Then, we get that 
\begin{align*}
    & \rho(M(.|x,u),\tau_{\sigma}(.|x,u))\leq \rho(M(.|x,u),M(.|\Bar{x}_{i},\Bar{u}_{i}))+ \rho(M(.|\Bar{x}_{i},\Bar{u}_{i}),a_{\sigma(i)})\\
    & \quad \quad \quad  \quad \quad+\rho(a_{\sigma(i)},\tau_{\sigma}(.|\Bar{x}_{i},\Bar{u}_{i}))+\rho(\tau_{\sigma}(.|\Bar{x}_{i},\Bar{u}_{i}),\tau_{\sigma}(.|x,u)) < \epsilon
\end{align*}

Thus, $ \sup_{(x,u)\in \mathbb{X}\times \mathbb{U}}\rho\big(M(.|x,u),\tau_{\sigma}(.|x,u)\big) \leq  \epsilon$ which implies that $\mathbb{M}$ is precompact and hence the closure of $\mathbb{M}$ under the uniform bounded Lipschitz metric is compact. 
\qed
{\begin{remark}
    We note that if Assumption \ref{CompactStateActionSpaces} does not hold but $\mathbb{X}\times \mathbb{U}$ is $\sigma$-compact, a similar proof to that of Theorem \ref{them exisistence of eps-partitions} can be applied to every compact subset of $\mathbb{X}\times\mathbb{U}$ to show that $\mathbb{M}$ is pre-compact under the topology of uniform weak convergence on compact subsets (see Definition \ref{unif weak on compact}).
\end{remark}}

\begin{remark}
    As demonstrated in the proof of Theorem \ref{them exisistence of eps-partitions}, uniform equicontinuity in Assumption \ref{equicontionuity} could be weakened to pointwise equicontinuity. 
\end{remark}
The following implies Assumption \ref{equicontionuity}.

\begin{proposition}
    Suppose $\mathbb{M}$ consists of the set of kernels which are realized by a measurable function $F(x,u,w)$ such that 
    \[X_{t+1}=F(X_{t},U_{t},W_{t})\] where $F\in \mathcal{A}$ and $\mathcal{A}$ is a family of equicontinuous functions over $x,u$ uniformly over $w$. Here, $W_{t}$ is some i.i.d noise process. Then, Assumption \ref{equicontionuity} holds.
\end{proposition}   
\noindent \textbf{Proof.} Let $\epsilon>0$, $(x,u)\in\mathbb{X}\times \mathbb{U}$. Let $\delta>0$ be such that for any $(x',u')$ such that $d\big((x,u),(x',u')\big)<\delta$: $d(F(x,u,w),F(x',u',w))\leq \epsilon$ for all $F\in\mathcal{A}$. Then, we get that for any $(x',u')$ such that $d\big((x,u),(x',u')\big)<\delta$ 
    \begin{eqnarray*}
\nonumber         \rho(M(.|x,u),M(.|x',u'))&=&\sup_{\|f\|_{BL}\leq 1}\mid\int f(F(x,u,w))P(dw)-\int f(F(x',u',w))P(dw)  \mid\\
 \nonumber       &\leq& \int d\big(F(x,u,w),F(x',u',w)\big)P(dw)\leq \epsilon
    \end{eqnarray*}
    \qed

\begin{proposition}
    Assumption \ref{assumption for proposition II} implies Assumption \ref{equicontionuity}.
\end{proposition}
\textbf{Proof.} Let $\epsilon>0$, $(x,u)\in \mathbb{X}\times \mathbb{U}$. Let $\delta>0$ be such that for any $(x',u')$ such that $d\big((x,u),(x',u')\big)<\delta$: $d(f^{M}(x,u,y),f^{M}(x',u',y))\leq \epsilon$ for all $M\in \mathbb{M}$. Then, we get that for any $(x',u')$ such that $d\big((x,u),(x',u')\big)<\delta$
\begin{align*}
     &\sup_{\|g\|_{BL}\leq 1}\bigg|\int g(y)M(dy|x,u)-\int g(y)M(dy|x',u')\bigg|\\
     &=\sup_{\|g\|_{BL}\leq 1}\bigg|\int g(y)f^{M}_{x,u}(y)\psi(dy)-\int g(y)f_{x',u'}^{M}(y)\psi(dy)\bigg| \\
     &\leq \int \sup_{M\in \mathbb{M}}\bigg|f^{M}_{x,u}(y)- f_{x',u'}^{M}(y)\bigg|\psi(dy)\leq \epsilon
\end{align*}
\qed

\begin{corollary}
    Assumption \ref{assumption for proposition II} and \ref{CompactStateActionSpaces} imply that the closure of $\mathbb{M}$ is compact. 
\end{corollary}

\subsection{Continuity and robustness of the optimal cost in the transition kernel} \label{Continuity of the optimal cost in the transition kernel}
In this subsection, we review recent results on the continuity of the optimal cost in the transition kernel \cite{kara2022robustness} which will be crucial for establishing near optimality of adaptive control policies.

\begin{assumption} \label{Strong ergodicity Ali}
    Let $\{\mathcal{T}_{n}\}_{n\in \mathbb{N}}\subseteq\mathbb{M}$ be such that for every stationary policy $\gamma\in \Gamma_{S}$, the transition kernels $\mathcal{T}_{n}$ and $\mathcal{T}$ lead to positive Harris recurrent chains. In particular, they lead to invariant probability measures $\pi^{n}_{\gamma}$ and $\pi_{\gamma}$ such that for any initial state $x\in \mathbb{X}$ we have that:
    \begin{eqnarray}
        \lim_{t\rightarrow\infty}\sup_{\gamma\in\Gamma_{S}}\|P_{t,x}^{\gamma,\mathcal{T}}(.)-\pi_{\gamma}(.)\|_{TV}&=&0 \label{equation 4.1}\\
        \lim_{t\rightarrow\infty}\sup_{n}\sup_{\gamma\in\Gamma_{S}}\|P_{t,x}^{\gamma,\mathcal{T}_{n}}(.)-\pi_{\gamma}^{n}(.)\|_{TV}&=&0
    \end{eqnarray}
    Where, $P_{t,x}^{\gamma,\mathcal{T}}$ denotes the probability induced on $X_{t}$ given that the system starts at $x$, the transition kernel is $\mathcal{T}$, and the policy is $\gamma$. Similarly, $P_{t,x}^{\gamma,\mathcal{T}_{n}}$ is the probability induced on $X_{t}$ given that the system starts at $x$, the transition kernel is $\mathcal{T}_{n}$, and the policy is $\gamma$.
\end{assumption}
Assumption \ref{Strong ergodicity Ali} is crucial because, under the additional Assumption that $c(x,u)$ is continuous and bounded, and weak continuity of the kernel in state and actions (Assumption \ref{AssumptionWFAverage}), it guarantees the existence of a measurable solution to the average cost optimality equation (ACOE) and in particular an optimal policy can be selected from $\Gamma_{S}$  (see Theorem 2.2 in \cite{hernandez2012adaptive}). 

Consider the Bellman operator: $\mathbb{T}: C_b(\mathbb{X}) \to C_b(\mathbb{X})$ with \[\mathbb{T}v(x) :=\inf_{u \in \mathbb{U}} \bigg( c(x,u)+\int_{\mathbb{X}}v(y)\mathcal{T}(dy|x,u)\bigg).\] Then, the ACOE is given by
\begin{equation} \label{ACOE}
    j^{*}+v^{*}(x)=\mathbb{T}v^{*}(x)
\end{equation}

If equation (\ref{ACOE}) has a solution, then under mild conditions \cite{survey}, $j^{*}$ represents the optimal cost and the optimal policy is the one that achieves the minimum in the right hand side of equation (\ref{ACOE}). For further conditions which guarantee the existence of a solution to (\ref{ACOE}), see \cite[Corollary 3.6]{hernandez2012adaptive}.
\begin{assumption} \label{Assumption for ACOE}
There exists $\beta \in [0,1)$ such that
\begin{eqnarray*}
    \sup_{M\in \mathbb{M}}\sup_{(x,u),(x,u')}\|M(.|x,u)-M(.|x',u')\|_{TV}&\leq& 2 \beta
\end{eqnarray*}
\end{assumption}
Assumption \ref{Assumption for ACOE} guarantees that for all $M\in \mathbb{M}$ a solution to the ACOE exists (corollary 3.6 in \cite{hernandez2012adaptive}). We note that Assumption \ref{Assumption for ACOE} implies Assumption \ref{Strong ergodicity Ali} by an application of Lemma 3.3 in \cite{hernandez2012adaptive}.

\begin{theorem} \cite[Theorem 2]{kara2022robustness}, \cite[Theorem 4.4]{kara2020robustness} \label{Theorem 2 Ali}
    Suppose Assumption \ref{Strong ergodicity Ali}, and Assumption \ref{AssumptionWFAverage} hold. Suppose that $\mathcal{T}_{n}\rightarrow\mathcal{T}$ weakly continuously. Additionally, suppose that a solution to the ACOE (\ref{ACOE}) holds for every $\mathcal{T}_{n}$. Then $J(\mathcal{T},X_0,\gamma_{n}^{*})\rightarrow J^{*}(\mathcal{T},X_{0})$ for any initial state $X_{0}\in \mathbb{X}$, where $\gamma_{n}^{*}$ is an optimal policy corresponding to the kernel $\mathcal{T}_{n}$ and $J(\mathcal{T},X_0,\gamma_{n}^{*})$ denotes the cost incurred when the policy $\gamma_{n}^{*}$ is applied. 
\end{theorem}




\subsection{Continuity of invariant probability measure in stationary Markov policies}

In this section, we study continuous dependence of invariant measures on control policy under Young topology. Such a continuous dependence will allow for asymptotic identifiability with an arbitrarily small perturbation loss. 

\textit{Definition}: Let $\mathcal{G}$ denote the set of \textit{invariant occupation measures}, also known as \textit{ergodic occupation measures}, which is defined by 
\begin{eqnarray*}
    &&\mathcal{G}:=\{\mu \in \mathcal{P}(\mathbb{X} \times \mathbb{U})| \mu(B\times \mathbb{U})=\int_{\mathbf{K}}\mu(dx,du)\mathcal{T}(B|x,u),\text{ for all }B\in\mathcal{B}(\mathbb{X})\}
\end{eqnarray*}
Consider the following assumptions regarding the true kernel $\mathcal{T}$.


\begin{theorem} \label{continuity in invariant measures}
    [Theorem 4.1 \cite{yuksel2023borkar}] Suppose that
    \begin{enumerate}
        \item[(a)] We have $\mathbb{X}=\mathbf{R}^{n}$ for some finite $n$, and for all $x\in\mathbf{R}^{n}$.
         \item[(b)] The kernel $\mathcal{T}(dy|x,u)$ is such that the family of conditional probability measures $\{\mathcal{T}(dy|x,u),x\in\mathbb{X},u\in\mathbb{U}\}$ admit densities $f_{x,u}(y)$ with respect to a reference probability measure $\psi \in \mathcal{P}(\mathbb{X})$, and all such densities are bounded and equicontinuous (over $x\in\mathbb{X}$, $u\in \mathbb{U}$).
        \item[(c)] For every stationary policy $\gamma \in \Gamma_{S}$ there exists a unique invariant probability measure.
        \item[(d)] Assumption \ref{AssumptionWFAverage} (b) holds and $\mathcal{G}$ is weakly compact (a sufficient condition being Assumption \ref{CompactStateActionSpaces}).
    \end{enumerate}
    Then, the invariant measure is weakly continuous on the control policy space $\Gamma_{S}$. Here, policies are endowed with the Young topology at input $\psi$ (See Definition \ref{definition young Top}).
\end{theorem}

{
\begin{remark} While the above is stated for $\mathbb{X} = \mathbb{R}^n$, the results also hold provided that $\mathbb{X}$ is standard Borel and $\sigma$-compact, and the other conditions remain applicable.
\end{remark}
}
As will be seen, the continuity result of Theorem \ref{continuity in invariant measures} leads to robustness in terms of the choice of the exploration policy needed to achieve learning. 

In the following two sections, we present two data driven identifiability and adaptive control results. The first one relies on performing Bayesian updates on the prior $\Theta$. The second identifiability result, on the other hand, relies on computing the empirical occupation measures.

\section{Bayesian Learning: Identifiability and Adaptive Control}

In this section, we consider Bayesian learning, identifiability and asymptotically optimal adaptive control. While in the literature the parametric or the finite cases have been studied under restrictive technical assumptions for the continuous space setting \cite{Kim2017TSfinite-case-average-cost,Adelman2024Thompsonsampling,Kim2019TS-continuous-case-average-cost}, our contribution is on the general non-parametric setup with relatively modest assumptions. The generality, however, prevents us from obtaining sample complexity bounds.

\subsection{Bayesian identifiability}

Given a prior probability measure $\Theta$ on $\mathbb{M}$, a stationary Markovian exploration policy $\gamma_{E}$, and an initial state distribution $X_{0}\sim \mu$ which is independent of $\Theta$, and any $T \in \mathbb{Z}_+$, one can define a unique probability measure on $\mathbb{M}\otimes_{k=0}^{T}\mathbb{X}\otimes_{k=0}^{T-1}\mathbb{U}$, which can be extended to the infinite horizon using Ionescu-Tulcea extension Theorem \cite{Ionescu-Tulcea}, given by 
\begin{align} \label{Bayesian update}
    P^{\Theta,\gamma_{E},\mu}(dM,dX_{[0,T]},dU_{[0,T-1]})=\int \prod_{k=0}^{T-1}(M(dX_{k+1}|X_{k},U_{k})\gamma_{E}(dU_{k}|X_{k}))\mu(dX_{0})\Theta(dM)
\end{align}

It thus follows that the (Bayesian) posterior $\Theta_{t}=P^{\Theta}(dM|I_{t})=P^{\Theta}(dM|X_{[0,t]},U_{[0,t-1]})$ is a well defined conditional probability measure. In our analysis below, we will not compute explicitly the conditional probability measure but our analysis will imply the asymptotic near recovery of the true kernel leading to near optimality. Notably, we will present conditions under which the controller can asymptotically, almost surely, and nearly identify the true transition kernel. For an explicit computation of the posterior above see \cite{Kim2017TSfinite-case-average-cost} for the case where $\mathbb{M}$ is finite, and \cite{Kim2019TS-continuous-case-average-cost} for the continuous case.

\begin{definition}
    Let $\epsilon>0$. We say that $\mathbb{M}_{\epsilon}$ is an $\epsilon$-net of $\mathbb{M}$ if and only if for some finite $N\in\mathbb{N}$, $\mathbb{M}_{\epsilon}=\{M_{1},...,M_{N}\}\subseteq \mathbb{M}$ is such that for all $M\in \mathbb{M}$ $\displaystyle \min_{M_j \in \mathbb{M}_{\epsilon}}\rho^{\mathrm{unif}}_{BL}(M,M_j) \leq \epsilon$.

    Thus $\mathbb{M}$ can be partitioned into bins $\{Q_{i}\}_{i=1}^{N}$ such that $Q_{i}=Q^{-1}(M_{i})$ where $Q$ is a quantizer function that maps every element of $M\in \mathbb{M}$ to the nearest element in $\mathbb{M}_{\epsilon}$. Thus $\bigcup_{i=1}^{N}Q_{i}=\mathbb{M}$ forms an $\epsilon$-partition of $\mathbb{M}$. It will be useful to make the bins disjoint by defining  $\Bar{Q}_{1}=Q_{1}$ and $\Bar{Q}_{j}=Q_{j}\bigcap_{k=1}^{j-1}\Bar{Q}_{k}^{c}$ for any $j\in \{2,...,N\}$. 
    
\end{definition}
The following is a crucial assumption which is assumed to hold throughout the paper. 
{
\begin{assumption} \label{prior charges bins}
 The true kernel $\mathcal{T}$ is in the topological support of $\Theta$. 
\end{assumption}
}


We note that Assumptions \ref{CompactStateActionSpaces} and \ref{equicontionuity} or Assumption \ref{compactness under unif BL} guarantee the existence of an $\epsilon$ net $\mathbb{M}_{\epsilon}$ for every $\epsilon>0$.

\begin{definition} \label{F definition}
    Let $\phi \in \mathcal{P}(\mathbb{X}\times\mathbb{U})$ be a probability measure on $\mathbb{X}\times\mathbb{U}$ and $\mathbf{F}$ be a countable family of continuous and bounded Lipschitz functions that can be used to characterize the weak topology \cite[Theorem 3.4.5]{ethier2009markov}. We say that two transition kernels $M$, $M'\in \mathbb{M}$ are $\phi$-equivalent and write $M\equiv M'$ if and only if for all $(x,u)$ in the topological support of $\phi$ (i.e., $(x,u)$ is such that for every open neighborhood $N_{(x,u)}$ of $(x,u)$, we have that $\phi(N_{(x,u)})>0$), the following holds for all $f\in \mathbf{F}$: 
    \begin{equation} \label{equivalence of kernels}
        \int f(s) M(ds|x,u)= \int f(s) M'(ds|x,u)
    \end{equation}
\end{definition}    

\begin{definition}
    Let $\mathcal{F}_{t}:=\sigma\big((X_{t},U_{t}),(X_{t+1},U_{t+1}),...\big)$. The tail sigma algebra is defined as: 
    $\displaystyle \mathcal{F}_{\infty}:=\bigcap_{t=0}^{\infty}\mathcal{F}_{t}$
\end{definition} 
{\begin{definition}
    A Markov chain $\{X_{t}\}_{t\in \mathbb{Z}^{+}}$ with state space $\mathbb{X}$ is said to be $\psi$-irreducible, for some $\psi\in\mathcal{P}(\mathbb{X})$, if and only if for every $A\in \mathcal{B}(\mathbb{X})$ such that $\psi(A)>0$, and any $x\in \mathbb{X}$ there exists $t\in \mathbb{N}$ such $P(X_{t}\in A|X_{0}=x)>0$. Furthermore,  the Markov chain is said to be periodic with period $d$ if there exists a partition  $X = \bigcup_{i=1}^{d} \mathbb{X}_{d}\cup D$ so
that $P(\mathbb{X}_{i+1}|x) = 1$ for all $x \in \mathbb{X}_{i}$ and $P(X_{1}|x) = 1$ for all $x \in \mathbb{X}_{d}$, with $\psi(D) = 0$. If no such $d > 1$ exists, the chain is aperiodic.
\end{definition}
\begin{definition}
    A Markov chain $\{X_{t}\}_{t\in \mathbb{Z}^{+}}$ with state space $\mathbb{X}$ is said to be $\phi$-recurrent, for some $\phi \in \mathcal{P}(\mathbb{X})$, if and only if for every $A\in \mathcal{B}(\mathbb{X})$ such that $\phi(A)>0$, and any $x\in A$ $P(\{X_{t}\in A \text{ for some }t\in \mathbb{N}\}|X_{0}=x)=1$.
\end{definition}}
\begin{definition}
    A policy $\gamma \in \Gamma_{S}$ is said to be \textbf{exhaustive} if and only if it induces an invariant probability measure $\Bar{\pi}$ on $\mathbb{X}\times \mathbb{U}$ such that for any open set $B\subseteq \mathbb{X}\times \mathbb{U}$ such that $B\not=\emptyset$: $\Bar{\pi}(B)>0$.
\end{definition}
\begin{theorem} \label{exhaustive policies are dense}
     Suppose Assumption \ref{Assumption for ACOE} holds for the true kernel. Moreover, suppose that the conditions of Theorem \ref{continuity in invariant measures} hold. If an exhaustive policy exists, then the set of all exhaustive policies is weakly dense in $\Gamma_{S}$. 
\end{theorem}

\textbf{Proof.} Let $\gamma\in \Gamma_{S}$ be an exhaustive policy and let $\gamma_{s}\in \Gamma_{S}$ be any stationary policy. Let $w\sim Ber(p_{n})$ and $\gamma_{n}=\begin{cases}
    \gamma  & \text{if } w=1\\
    \gamma_{s}  & \text{if } w=0
\end{cases}$. 

\begin{eqnarray*}
 &&\sup_{\|f(x,u)\|_{\infty}\leq 1}\bigg| E^{\gamma_{n}}[f(x,u)]- E^{\gamma_{s}}[f(x,u)]  \bigg| \leq (1-p_{n})\sup_{\|f\|_{\infty}\leq 1}\bigg| E^{\gamma_{s}}[f(x,u)|w=0]-E^{\gamma_{s}}[f(x,u)]\bigg|  \\
\nonumber  \\  && \qquad +  p_{n}\sup_{\|f\|_{\infty}\leq 1}\bigg| E^{\gamma_{n}}[f(x,u)|w=1]-E^{\gamma_{s}}[f(x,u)]\bigg| \leq 2p_{n}.
\end{eqnarray*}
Thus, as $p_{n}\rightarrow 0$, $\sup_{\|f(x,u)\|_{\infty}\leq 1}\bigg| E^{\gamma_{n}}[f(x,u)]- E^{\gamma_{s}}[f(x,u)]  \bigg|\rightarrow 0$ (pointwise in $x$), and hence for any continuous and bounded function $f$ \[\bigg| \int f(x,u)\gamma_{n}(du|x)\psi(dx)-\int f(x,u)\gamma(du|x)\psi(dx)\bigg| \rightarrow 0 \] Hence,  $\gamma_{n}\rightarrow \gamma_{s}$ weakly at input measure $\psi$ (recall Definition \ref{definition young Top}). 

Thus, it follows from Theorem \ref{continuity in invariant measures} that $\Bar{\pi}_{n}\rightarrow \Bar{\pi}_{s}$ weakly. Where, $\Bar{\pi}_{n}$ and $\Bar{\pi}_{s}$ are the invariant probability measures induced by $\gamma_{n}$ and $\gamma_{s}$ respectively. Next, we will show that for any open set $B\in \mathbb{X}\times \mathbb{U}$ such that $B\not=\emptyset$ $\Bar{\pi}_{n}(B)>0$. Consider an open set $B\in \mathbb{X}\times \mathbb{U}$ such that $B\not=\emptyset$. Since $\gamma$ is an exhaustive policy, it follows that its induced invariant probability measure $\Bar{\pi}$ satisfies $\Bar{\pi}(B)>0$. By Assumption \ref{Assumption for ACOE} it follows that $\sup_{x}\|P^{t,\gamma}_{x}-\Bar{\pi}\|_{TV}\rightarrow 0$
as $t\rightarrow\infty$. In particular, this implies that $P^{t,\gamma}_{x}(B)\rightarrow \Bar{\pi}(B)>0$ uniformly over $x$. Let $\Bar{\pi}(B)>\epsilon>0$. There exists $T>0$ such that for all $t\geq T$ and for all $x$ $P^{t,\gamma}_{x}(B)>\epsilon$. Hence, we get that 
\begin{eqnarray*}
    \Bar{\pi}_{n}(B)= \int \Bar{\pi}_n(dx) P^{T,\gamma_n}_x(B) \geq (p_n)^T \inf_{x\in \mathbb{X}}  P_x^{T,\gamma}(B) \geq  (p_n)^T \epsilon > 0.
\end{eqnarray*}
\qed

{Let $M$ be random variable which denotes the true kernel. Before we proceed with the main identifiability theorem we present a useful lemma which holds when $\mathbb{X}$ and $\mathbb{U}$ are $\sigma$-compact.
\begin{lemma} \label{M is tail measurable}
    Suppose Assumption \ref{equicontionuity} holds and, in addition, either Assumption \ref{CompactStateActionSpaces} or Assumption \ref{compactness under unif BL} holds. Additionally, suppose $\gamma$ is an exhaustive exploration policy, i.e., an exhaustive policy applied by the controller to learn the transition kernel, which for every $m\in\mathbb{M}$ induces a $\psi$-irreducible, $\phi$-recurrent, and aperiodic MC with invariant probability measure $\Bar{\pi}$. Then, $M$ is measurable with respect to $\mathcal{F}_{\infty}$ and hence measurable with respect to $I_{\infty}:= \sigma(X_{0},U_{0},X_{1},U_{1},...)$ $P^{\Theta}$ almost surely where $P^{\Theta}$ denotes the joint measure induced on the entire state action process $\{(X_{t},U_{t})\}_{t=0}^{\infty}$ and the space of kernels $\mathbb{M}$ under the exploration policy and given the prior $\Theta$.
\end{lemma}

}

\textbf{Proof}. See Appendix.

The following theorem also holds when $\mathbb{X}$ and $\mathbb{U}$ are $\sigma$-compact.
\begin{theorem} \label{martingale convergent}
 Suppose Assumption \ref{equicontionuity} holds and, in addition, either Assumption \ref{CompactStateActionSpaces} or Assumption \ref{compactness under unif BL} holds. Additionally, suppose $\gamma$ is an exhaustive exploration policy, i.e., an exhaustive policy applied by the controller to learn the transition kernel, which for every $m\in\mathbb{M}$ induces a $\psi$-irreducible, $\phi$-recurrent, and aperiodic MC with invariant probability measure $\Bar{\pi}$. Then, for every $\epsilon>0$, with $\displaystyle M_{t}= \arg\max_{\{M_{j}\in \mathbb{M}_{\epsilon}\}} P^{\Theta}(M\in \Bar{Q}_{j}|I_{t})$, where $P^{\Theta}(M\in \Bar{Q}_{j}|I_{t})$ denotes the posterior distribution on $\mathbb{M}$ given that the prior is $\Theta$, satisfies for all $(x,u) \in \mathbb{X} \times \mathbb{U}$: $\rho\big(M_{t}(.|x,u),\mathcal{T}(.|x,u)\big)\leq \epsilon$ as $t\rightarrow \infty$ almost surely.
\end{theorem}
\begin{remark}
    It is important to note that the measures $\phi$, $\psi$, and $\Bar{\pi}$ in the statement of Lemma \ref{M is tail measurable} and Theorem \ref{martingale convergent} need not be the same for all $m\in \mathbb{M}$.
\end{remark}

\textbf{Proof}. \textbf{Step 1}: [\textit{Convergence of the posterior}] Let $\epsilon>0$ and $M_{j}\in \mathbb{M}_{\epsilon}$. Let $Q_{j}=Q^{-1}(M_{j})$. Define $\Bar{Q}_{1}=Q_{1}$ and $\Bar{Q}_{j}=Q_{j}\bigcap_{k=1}^{j-1}\Bar{Q}_{k}^{c}$ for any $j\in \{2,...,N\}$. Consider $P^{\Theta}(\{M\in \Bar{Q}_{j}\}|I_{t})$. We have that $P^{\Theta}(\{M\in \Bar{Q}_{j}\}|I_{t})$ is a bounded martingale with respect to the filtration $I_{t}$, hence we get that $P^{\Theta}(\{M\in \Bar{Q}_{j}\}|I_{t})$ converges a.s to some r.v $Z^{\ast}_{\epsilon}(\Bar{Q}_{j})$.  

\textbf{Step 2}: [\textit{The event that the kernel is in a bin is in the tail sigma algebra}] By Lemma \ref{M is tail measurable}, we have that $\{M\in Q_{j}\}\in \mathcal{F}_{\infty}$ where $\mathcal{F}_{\infty}$ denotes the tail sigma algebra. This, in turn, entails that $\{M\in \Bar{Q}_{j}\}\in \mathcal{F}_{\infty}$ by the definition of $\Bar{Q}_{j}$. 

\textbf{Step 3}: { [Consistency of the posterior with the true transition kernel]
We have that $Z^{*}_{\epsilon}(\Bar{Q}_{j})=P^{\Theta}(\Bar{Q}_{j}|I_{\infty})\in \{0,1\}$ where $I_{\infty}=\sigma(\{(X_{t},U_{t})\}_{t=0}^{\infty})$. Through iterated expectations we get that

\begin{eqnarray*}
  P^{\Theta}(\Bar{Q}_{j}|I_{\infty})&=&E^{P^{\Theta}}[\chi_{\{\Bar{Q}_{j}\}}|I_{\infty}]= E[E^{P^{\Theta}}[\chi_{\{\Bar{Q}_{j}\}}|I_{\infty},M]| I_{\infty}]\\
  &=&E[P^{\Theta}(\Bar{Q}_{j}|I_{\infty},M)| I_{\infty}]= E[\delta_{M}(\Bar{Q}_{j})|I_{\infty}]= \delta_{M}(\Bar{Q}_{j}).
\end{eqnarray*}
 Where the fourth equality follows from the fact that $P^{\Theta}(\Bar{Q}_{j}|M)=\delta_{M}(\Bar{Q}_{j})$ and the last equality follows from the fact that $M$ is measurable with respect to $I_{\infty}$ (by Lemma \ref{M is tail measurable}).

 Hence, almost surely for every $\epsilon>0$, there exists $T$ such that for all $t\geq T$: $M_{t}\in \Bar{Q}_{j}^{-1}(M_{j})$. Where $M_{j} \in \mathbb{M}_{\epsilon}$ is such that for all $x,u$: $\rho\big(M_{j}(.|x,u),\mathcal{T}(.|x,u)\big)\leq \epsilon$.} 
\qed 

\subsection{Bayesian adaptive control}

\subsubsection{Alternating between exploration and exploitation}

Here we describe a strategy where the controller alternates between exploration periods and periods of exploitation. In particular, let 
 $T',T_{l}(T'),T_{a}(T')\in \mathrm{Z}^{+}$, where $T'=T_{a}+T_{l}$. The controller uses an exhaustive exploration policy which induces an irreducible, recurrent, and aperiodic MC, in the periods $t=kT',...,t=kT'+T_{l}-1$. Here, $k\in \mathrm{Z}^{+}$. In every subsequent exploration period, the controller attempts to improve on their kernel estimate but only considers the kernels which are within an $\epsilon$-distance of the initial estimate. This guarantees that all subsequent estimates are within a distance $2\epsilon$ of the true kernel. At time $t=kT'+T_{l}$ the controller estimates the maximum likelihood transition kernel using $\{X_{[KT',kT'+T'_{l}]},U_{[KT',kT'+T_{l}-1]}\}$ and applies the corresponding optimal stationary Markov policy in the period $t=kT'+T_{l},...,(k+1)T'-1$. We denote such adaptive policies by $\gamma(\epsilon,T')$.

\begin{theorem} \label{adavptive control theorem}
    Suppose the following holds
    \begin{itemize}
        \item[(i)] Assumption \ref{all is WF}.
        \item[(ii)] Assumption \ref{equicontionuity} and Assumption \ref{CompactStateActionSpaces}, or Assumption \ref{compactness under unif BL} hold.
        \item[(iii)] Assumption \ref{Assumption for ACOE} holds.
        \item[(iv)] $\displaystyle \lim_{T'\rightarrow \infty}\frac{T_{l}}{T'}=0$ and $\displaystyle \lim_{T'\rightarrow \infty} T_{l}=\infty$
    \end{itemize}
   Then, the exploitation-exploration policy described above is near optimal for sufficiently large $T'$.
\end{theorem}
\textbf{Proof.} 
 Step 1: Let $\epsilon>0$ and $\epsilon'>0$. Let $\nu$ denote the probability measure induced on the entire state process $\{(X_{0},U_{0}),(X_{1},U_{1}),...\}$. Note that the existence of such a measure is guaranteed by Ionescu-Tulcea's Theorem \cite{Ionescu-Tulcea}. By Theorem \ref{martingale convergent} we know that the maximum likelihood estimator converges $\nu$ almost surely. Hence, by Egorov's Theorem \cite{Egorovtheorem1911} for any $\delta>0$ there exists a set $B$ such that $\nu(B^{C})<\delta$ and on $B$ the convergence of the estimator is uniform. Hence, on the set $B$, we can choose $T_{l}$ large enough so that the initial exploration period is sufficient to learn the true transition kernel to within some error tolerance $\epsilon>0$ and the value of $T_{l}$ is independent of the sample path $\omega$. In particular, because the cost functional $c(x,u)$ is bounded, we get that 
 \begin{align*}
      \nonumber  J(\gamma)&=\limsup_{T\rightarrow \infty} E^{\gamma}_{X_{0}}[\frac{1}{T}\sum_{t=0}^{T-1}c(X_{t},U_{t})]\\
      &\leq \limsup_{T\rightarrow \infty}\bigg( E^{\gamma}_{X_{0}}[\frac{1}{T}\sum_{t=0}^{T-1}c(X_{t},U_{t})|B]+ E^{\gamma}_{X_{0}}[\frac{1}{T}\sum_{t=0}^{T-1}c(X_{t},U_{t})|B^{c}]\bigg)\\
      &\leq \limsup_{T\rightarrow \infty} E^{\gamma}_{X_{0}}[\frac{1}{T}\sum_{t=0}^{T-1}c(X_{t},U_{t})|B]+ \delta \|c\|_{\infty}
 \end{align*}
 
Thus, since we are only interested in near optimal policies, one can without loss of generality assume that the convergence in Theorem \ref{martingale convergent} is uniform over all sample paths. Thus one can construct an exploitation-exploration policy $\gamma(\epsilon',T')$ such that in the initial exploration phase the kernel is identified up to some $\epsilon'$.

Step 2: We have that 
\begin{eqnarray*}
  \nonumber  J(\gamma)&=&\limsup_{T\rightarrow \infty} E^{\gamma}_{X_{0}}[\frac{1}{T}\sum_{t=0}^{T-1}c(X_{t},U_{t})]=\limsup_{n\rightarrow \infty} E^{\gamma}_{X_{0}}[\frac{1}{nT'}\sum_{t=0}^{nT'-1}c(X_{t},U_{t})]\\
 \nonumber   &\leq& \limsup_{n\rightarrow \infty} \frac{1}{nT'} \sum_{k=0}^{n-1}\big[\sum_{t=kT'}^{kT'+T_{l}-1}E[c(X_{t},U_{t})]+\sum_{t=kT'+T_{l}}^{(k+1)T'-1} E[c(X_{t},U_{t})]\big]\\
 \nonumber   &\leq& \limsup_{n\rightarrow \infty} \frac{1}{nT'} \sum_{k=0}^{n-1}\sum_{t=kT'}^{kT'+T_{l}-1}E[c(X_{t},U_{t})]+\limsup_{n\rightarrow \infty} \frac{1}{nT'} \sum_{k=0}^{n-1}\sum_{t=kT'+T_{l}}^{(k+1)T'-1} E[c(X_{t},U_{t})]\\
 \nonumber   &\leq& \|c\|_{\infty}\frac{T_{l}}{T'}+\limsup_{n\rightarrow \infty} \frac{T_{a}}{nT'} \sum_{k=0}^{n-1}\sum_{t=kT'+T_{l}}^{(k+1)T'-1}\frac{1}{T_{a}} E[c(X_{t},U_{t})]
\end{eqnarray*}
Since $\displaystyle \lim_{T'\rightarrow \infty}\frac{T_{l}}{T'}=0$, the first term in the bound above can be made arbitrarily small. 

Step 3: Next, we bound the second term. Let $k\in \{0,...,n-1\}$. 


By the continuity result of Theorem \ref{Theorem 2 Ali} there exists $\delta>0$ such that $\rho(M(.|x,u),\mathcal{T}(.|x,u))\leq \delta$ for all $x,u$ implies that $|J(\gamma(M),\mathcal{T})-J(\gamma(\mathcal{T}),\mathcal{T})|\leq \frac{\epsilon}{4}$. Where, $\gamma(M)$ and $\gamma(\mathcal{T})$ denote the optimal policies designed for the transition kernels $M$ and $\mathcal{T}$ respectively. By \cite[Lemma 1]{kara2022robustness} there exists some $\Bar{T}>0$ such that for all $\gamma\in\Gamma_{S}$, we have that $T>\Bar{T}$ implies that $|J^{T}(\gamma)-J(\gamma)|\leq \frac{\epsilon}{4}$. Where $J^{T}(\gamma)$ denotes the finite horizon cost criteria: $J^{T}(\gamma)= \frac{1}{T}E^{\gamma}_{X_{0}}[\sum_{t=0}^{T-1}c(X_{t},U_{t})]$. Since, $T_{a}(T')\rightarrow \infty$ as $T'\rightarrow \infty$, there exists $T^{*}$ such that for all $T'\geq T^{*}$, we have that $T_{a}>\Bar{T}$.


Step 4: By Theorem \ref{martingale convergent}, one can construct a sequence of exploitation-exploration policies such $\{\gamma_{j}(T'_{j}, \delta_{j})\}_{j\in \mathbb{N}}$ such that $\delta_{j}\leq \delta$ and $T'_{j}\geq T^{*}$ for all $j$. Let $M_{k}$ denote the estimated kernel at time $t=kT'+T_{l}$ Thus for any element of the sequence $\{\gamma_{j}(T'_{j}, \delta_{j})\}_{j\in \mathbb{N}}$, we get that
\begin{eqnarray*}
\nonumber \sum_{t=kT'+T_{l}}^{(k+1)T'-1}\frac{1}{T_{a}} E[c(X_{t},U_{t})]&=& J^{T_{a}}(\gamma_{M_{k}}) \leq J(\gamma_{M_{k}})+\frac{\epsilon}{4}\leq J(\gamma_{\mathcal{T}})+\frac{\epsilon}{4}+\frac{\epsilon}{4}= J(\gamma_{\mathcal{T}})+\frac{\epsilon}{2}
\end{eqnarray*}
Step 5: Thus, we get for all $j\in \mathbb{N}$
\begin{eqnarray*}
\nonumber    J(\gamma_{j})&\leq&\|c\|_{\infty}\frac{T_{l}(T'_{j})}{T'_{j}}+\limsup_{n\rightarrow \infty} \frac{T_{a}(T'_{j})}{nT'_{j}} \sum_{k=0}^{n-1}\sum_{t=kT'+T_{l}}^{(k+1)T'-1}\frac{1}{T_{a}(T'_{j})}E[c(X_{t},U_{t})]\\
\nonumber    &\leq&\|c\|_{\infty}\frac{T_{l}(T'_{j})}{T'_{j}}+\limsup_{n\rightarrow \infty} \frac{T_{a}(T'_{j})}{nT'_{j}} \sum_{k=0}^{n-1}J(\gamma_{\mathcal{T}})+\frac{\epsilon}{2}\\
 \nonumber   &=&\|c\|_{\infty}\frac{T_{l}(T'_{j})}{T'_{j}}+\frac{T_{a}(T'_{j})}{T'_{j}}(J(\gamma_{\mathcal{T}})+\frac{\epsilon}{2})
\end{eqnarray*}
Now, by taking the limit as $j\rightarrow \infty$ we get that $ \lim_{j\rightarrow \infty}J(\gamma_{j})\leq J(\gamma_{\mathcal{T}})+\frac{\epsilon}{2}$
Thus, there exists $N\in \mathbb{N}$ such that for all $j\geq N$
$J(\gamma_{j})\leq J(\gamma_{\mathcal{T}})+\epsilon$
 \qed
 \begin{remark}
     It is easy to check that Theorem \ref{adavptive control theorem} still holds if the exploitation-exploration policies are defined in a manner such that the kernel estimates depend not only on the data acquired during learning periods but also the data acquired during the exploitation periods.
 \end{remark}

\subsubsection{Simultaneous exploration and exploitation}

 The objective of this section is to establish the near-optimality of policies that simultaneously learn the true kernel and exploit the knowledge on the kernel at every time stage. Throughout this section we will assume that the space of kernels $\mathbb{M}$ is finite and show that such policies are optimal. We note however that, as studied earlier, for the case where $\mathbb{M}$ is not finite, under the conditions reported, one can first obtain an $\epsilon$-partition, and then apply a similar process to achieve a near optimal cost. \\

\textbf{Description of simultaneous exploration and exploitation policies}

 Consider an exhaustive exploration policy $\gamma_{e}$ which induces an irreducible and aperiodic MC. Let $\gamma^{*}_{0}$ denote the optimal policy corresponding to the maximum likelihood estimator at time zero and suppose it is applied at $t=0$. Let $k$ be a random variable denoting the number of times that maximum likelihood estimator changes. Let $\{v_{t}\}_{t\in \mathbb{N}}$ be a sequence of independent random variables such that $v_{t}\sim Ber(\frac{1}{(1+k)^{2}})$. For the first exploration and exploitation phase: $t=1,...,t=T_{1}$ the controller applies the policy $\gamma_{t}=\begin{cases}
    \gamma_{e}  & \text{if } v_{t}=1\\
    \gamma_{0}^{*}  & \text{ otherwise } 
\end{cases}$. Here, $\displaystyle T_{1}=\inf\{t\in \mathbb{Z}^{+}| \exists M \in \mathbb{M}\text{ such that }P(M_{t}\equiv M|I_{t})\geq \frac{1}{2}\}$. As will be seen in Lemma \ref{useful lemma}, it  follows from Theorem \ref{martingale convergent} that $P(T_{1}<\infty)=1$. Similarly one can define a sequence of stopping times $\{T_{i}\}_{i\in \mathbb{N}}$ such that $\displaystyle T_{i}=\inf\{t>T_{i-1}| \exists M \in \mathbb{M}\text{ such that }P(M_{t}\equiv M|I_{t})\geq 1-\frac{1}{2^{i}}\}$. In each phase $i$ corresponding to $t=T_{i-1}+1,...,T_{i}$, the controller applies the policy $\gamma_{t}=\begin{cases}
    \gamma_{e}  & \text{if } v_{t}=1\\
    \gamma_{i-1}^{*}  & \text{ otherwise }
\end{cases}$ where $\gamma_{i-1}^{*}$ is the optimal policy corresponding to the maximum likelihood estimator at time $T_{i-1}$. 

\begin{theorem} \label{sim expl and expl thm}
    Suppose Assumption \ref{all is WF} holds and Assumption \ref{assumption for quantizing actions} holds for some positive measure $\zeta$ with $\zeta(\mathbb{X})>0$. Suppose that $\gamma_{e}$ induces a $\phi$-recurrent MC such that $\phi\ll\zeta$. Then, the simultaneous exploration and exploitation policy described above is optimal. 
\end{theorem}

\begin{lemma} \label{useful lemma}
Suppose the stationary policy $\gamma_{e}$ induces a $\psi$-irreducible, $\phi$-recurrent, and aperiodic MC. Additionally, suppose Assumption \ref{assumption for quantizing actions} for some positive measure $\zeta$ with $\zeta(\mathbb{X})>0$. Then, $\gamma_{t}=\begin{cases}
    \gamma_{e}  & \text{if } v_{t}=1\\
    \gamma_{i-1}^{*}  & \text{ otherwise }
\end{cases}$ also induces a $\psi$-irreducible, $\phi$-recurrent, and aperiodic MC.    
\end{lemma}
\textbf{Proof.} It follows from the proof of Theorem \ref{exhaustive policies are dense} that the MC induced by $\gamma_{t}$ is also $\psi$-irreducible and aperiodic. Moreover, since the kernel is minorized it also follows that $\gamma_{t}$ induces a $\phi$-recurrent MC.


\qed

\textbf{Proof of Theorem \ref{sim expl and expl thm}}
By Lemma \ref{useful lemma} and Theorem \ref{exhaustive policies are dense}, it follows that during every phase of exploration and exploitation the policy applied satisfies the conditions of Theorem \ref{martingale convergent}. Thus all the stopping times $T_{i}$ are finite almost surely. Hence, the maximum likelihood estimator will converge almost surely to the true kernel after some finite number of phases. Let $T_{i_{n}}$ be such that for all $t\geq T_{i_{n}}$ $M_{t}=\mathcal{T}$. We have that $\displaystyle \sum_{t=T_{i_{n}}}^{\infty}P(v_{t}=1)=\sum_{t=T_{i_{n}}}^{\infty}\frac{1}{(1+k)^{2}}<\infty$ . Hence, by Borel-Cantelli's Lemma, it follows that there exists some $T'> T_{i_{n}}$ such that for all $t\geq T:$ $v_{t}=0$. Hence, after some finite time $T$, the controller always applies the optimal policy. Let $\gamma_{SE}$ denote the simultaneous exploration and exploitation policy. We have that
\begin{eqnarray*}
    &&J(\gamma_{SE})=\limsup_{T\rightarrow \infty} E^{\gamma_{SE}}[\frac{1}{T}\sum_{t=0}^{T-1} c(X_{t},U_{t})]\\
    &&=\limsup_{T\rightarrow \infty} \sum_{l=0}^{\infty}E^{\gamma_{SE}}[\frac{1}{T}\sum_{t=0}^{T-1} c(X_{t},U_{t})\mid T'=l]P(T'=l)\\
    &&\leq  \sum_{l=0}^{\infty}\limsup_{T\rightarrow \infty}\frac{1}{T}E^{\gamma_{SE}}[\sum_{t=0}^{T-1} c(X_{t},U_{t})\mid T'=l]P(T'=l)\\
    &&\leq \sum_{l=0}^{\infty}\limsup_{T\rightarrow \infty}\frac{1}{T}\sum_{t=0}^{l}E^{\gamma_{SE}}[ c(X_{t},U_{t})\mid T'=l]P(T'=l)\\
    && \quad \quad \quad \quad \quad +\limsup_{T\rightarrow \infty}\frac{1}{T}\sum_{t=l+1}^{T-1}E^{\gamma_{SE}}[ c(X_{t},U_{t})\mid T'=l]P(T'=l)\\
    &&= \sum_{l=0}^{\infty}\limsup_{T\rightarrow \infty}\frac{1}{T}\sum_{t=l+1}^{T-1}E^{\gamma^{*}}[ c(X_{t},U_{t})\mid T'=l]P(T'=l)\\
    &&\leq \sum_{l=0}^{\infty}\limsup_{T\rightarrow \infty}\frac{1}{T}\sum_{t=l+1}^{T+l}E^{\gamma^{*}}[ c(X_{t},U_{t})\mid T'=l]P(T'=l)\\
    &&=\sum_{l=0}^{\infty} J_{l}^{T'}(\gamma^{*})P(T'=l)    =J(\gamma^{*})
\end{eqnarray*}
Here the first inequality follows from a generalized version of Fatou's lemma (Theorem 2.2 \cite{fatou-lemma}). $J_{l}^{T'}(\gamma^{*})$ denotes the expected average cost when the optimal policy $\gamma^{*}$ is applied and the initial state is sampled from $P(dX_{l+1}|T'=l)$. By the Assumption that $\mathcal{T}(.|x,u)\geq \zeta(.)$ for all $x,u$ for some positive measure $\zeta$ with $\zeta(\mathbb{X})>0$, it follows that the expected average cost is independent of the initial distribution (Lemma 3.3 \cite{hernandez2012adaptive}). Hence, one is able to obtain the last equality which then entails that $J(\gamma_{SE})=J(\gamma^{*})$.
\qed

 A significant challenge in the implementation of the adaptive policy presented in this section, when $\mathbb{M}$ is a continuous space, consists in the computation of the belief update on the space of transition kernel. The next section aims to circumvent this issue through a simple to implement empirical adaptive control algorithm which leads to near optimal solutions.  
   
\section{Empirical Learning: Identifiability and Adaptive Control}

As noted earlier, there have been several studies on empirical learning including those which learn noise from data or which obtain a quantized Markov approximation. Accordingly, this section builds on prior work towards an adaptive learning and control algorithm. Reliance on empirical occupation measures for kernel estimation is commonly used in the case where the noise is assumed to be observable (see for example \cite{yichen2024,Zurek-plugin-approach}). See also \cite[Section 4]{yichen2024} for an algorithm that learns an approximate model, assuming Lipschitz regularity of the kernel, through quantization and see \cite[Section 5.3]{ky2026qaverage} which establishes connections between empirical model learning and quantized Q-learning.

\subsection{The quantization procedure}
 Throughout this section, we assume that Assumption \ref{CompactStateActionSpaces} holds. Here, we describe a quantization procedure which was first introduced in \cite{SaYuLi15c}. Let $n \in \mathbb{N}$ and $\mathbb{U}^{m}=\{u_{1},u_{2},...,u_{m}\}$ be such that for all $u\in\mathbb{U}$ there exists $i\in \{1,...,m\}$ such that $d_{\mathbb{U}}(u,u_{i})<\frac{1}{n}$ where $d_{\mathbb{U}}$ is the metric on the space $\mathbb{U}$. Consider a MDP, which we denote $MDP_{n}=(\mathbb{X},\mathbb{U}^{m},\mathcal{T},c)$ such that the action space $\mathbb{U}$ is replaced with $\mathbb{U}^{m}$. We denote the optimal cost for $MDP_{n}$ by $J^{*}_{n}$. 
 \begin{assumption} \label{assumption for quantizing actions} The kernel $\mathcal{T}$ is minorized, i.e., there exists a measure $\zeta\in\mathcal{P}(\mathbb{X})$ such that for any $(x,u)\in\mathbb{X}\times\mathbb{U}$ and for any $B\in\mathcal{B}(\mathbb{X})$: $\mathcal{T}(B|x,u)\geq \zeta(B)$.  
 
 \end{assumption}
\begin{theorem} \label{action quantization theorem}[Theorem 3.16 in \cite{SaLiYuSpringer}]
Suppose Assumption \ref{assumption for quantizing actions} holds. Then, for any $x\in \mathbb{X}$: $|J^{*}_{n}(x)-J^{*}(x)|\rightarrow 0$ as $n\rightarrow \infty$.  
\end{theorem}
Next, we will describe a procedure for quantizing the state space. Assume $\mathbb{U}$ is finite or has been discretized according to the quantization procedure described above. Let $n \in \mathbb{N}$ and $\mathbb{X}^{m}=\{x_{1},x_{2},...,x_{m}\}$ be such that for all $x\in\mathbb{X}$ there exists $i\in \{1,...,m\}$ such that $d_{\mathbb{X}}(x,x_{i})<\frac{1}{n}$ where $d_{\mathbb{X}}$ is the metric on the space $\mathbb{X}$. Let $\nu \in \mathcal{P}(\mathbb{X})$ be a measure such that $\forall i \in\{1,...,n\}$ $\nu(B_{i})>0$ where $B_{i}=\{x\in \mathbb{X}| Q(x)=x_{i}\}$. One choice for $\nu$ that is often useful is $\displaystyle \nu(.)=\frac{1}{n}\sum_{i=1}^{n}\delta_{x_{i}}(.)$ where $\delta_{x_{i}}$ is the Dirac measure concentrated at $x_{i}$. Here, $Q$ is the quantizer that maps every element $x$ to the nearest element in $\mathbb{X}^{m}$. If ties occur, it is assumed that they are broken by whichever way makes the function $Q$ measurable. Consider the following transition kernel: 
\begin{equation} \label{P_{n}}
    P_{n}(x_{m}|x_{i},u')=\int_{B_{i}}\int_{\mathbb{X}} \chi_{\{Q(s)=x_{m}\}}\mathcal{T}(ds|w,u')\nu_{i}(dw)
\end{equation} Here $\nu_{i}$ denotes the restriction of the measure $\nu$ to the set $B_{i}$. Consider the MDP given by $MDP_{n}=(\mathbb{X}^{m},\mathbb{U},P_{n},c)$.  We denote by $J^{*}_{n}$ the expected cost accrued when the optimal policy for $MDP_{n}$ is extended into a policy for the original MDP which is constant over each bin.

\begin{theorem} [Theorem 3.7 \cite{ky2026qaverage}] \label{theorem sate space quantization}
    Under Assumption \ref{assumption for quantizing actions}, $|J^{*}_{n}(x)-J^{*}(x)|\rightarrow 0$ as $n\rightarrow \infty$.   
\end{theorem}


\subsection{The empirical learning algorithm}

\begin{algorithm}[h]
\caption{Algorithm I}
\label{alg:alg1}
\begin{algorithmic}[1]

\State \textbf{Quantization:}
\State For some fixed quantization step $n$ simultaneously quantize the state and action space.
\State Denote the resulting MDP by $MDP_{n}$ with kernel $P_{n}$.
\State Denote the resulting state and action spaces by 
$\mathbb{X}^{m}=\{x_{1},...,x_{m}\}$ and 
$\mathbb{U}^{m}=\{u_{1},...,u_{m}\}$.
\State \textbf{Initialization:}
\State Set $t=0$.
\State Let $M_{0}$ be any arbitrary kernel estimate.


\While{$t \leq T$}
\State Select an action uniformly at random from $\mathbb{U}^{m}$: $U_{t}\sim \mathrm{Unif}
(\mathbb{U}^{m})$
\State $t := t+1$
\EndWhile
\State Compute the kernel 
        $M_{T}:\mathbb{X}^{m}\times \mathbb{U}^{m}\rightarrow\mathcal{P}(\mathbb{X}^{m})$
        \State For any $(i,j,k)\in\{1,...,m\}^{3}$ define $ M_{T}(x_{k}\mid x_{i},u_{j})=$
        \[         \frac{1}{\sum_{m=0}^{T-1}\chi_{\{X_{m}=x_{i},U_{m}=U_{j}\}}}\sum_{l=1}^{T}\chi_{\{X_{l}=x_{k},U_{l-1}=x_{i},U_{l-1}=u_{j}\}}.
        \]
\State \Return{ $\gamma_{T}$ the optimal policy corresponding to $M_{T}$}
\end{algorithmic}
\end{algorithm}

In view of Theorem \ref{action quantization theorem}, we will assume that the action set $\mathbb{U}$ is finite in the following. We will consider Algorithm \ref{alg:alg1}.

\begin{theorem} \label{computational learning theorem I}
    Consider Algorithm \ref{alg:alg1}. Suppose that Assumption \ref{assumption for quantizing actions} holds. Additionally, suppose that $\zeta\in\mathcal{P}(\mathbb{X})$ satisfies $\zeta(B_{i})>0$ for each bin $B_{i}$. Then, as $T\rightarrow \infty$, we have that $M_{T} \rightarrow P_{n}$ weakly almost surely for some transition kernel of the form (\ref{P_{n}}). Moreover, $J(\gamma_{T}^{*}) \rightarrow J^{*}(\mathcal{T})$ almost surely where $J^{*}(\mathcal{T})$ is the optimal performance when the true kernel is known and $\gamma^{*}_{T}(x):=\sum_{X_{i}\in\mathbb{X}^{m}} \chi_{\{Q(x)=x_{i}\}} \gamma_{T}(x_{i})$.
\end{theorem}

\textbf{Proof.} Throughout the proof we will assume without any loss of generality that after the quantization process the quantized states and actions have the same cardinality. \\

\noindent Step 1: By Assumption \ref{assumption for quantizing actions}, we have that for any $x,x' \in \mathbb{X}^{m}$ and $u'\in \mathbb{U}^{m}$
\begin{eqnarray} \label{minorization-partial-inequality}
  \nonumber P_{n}(x|x',u')&=&\int_{B_{i}}\int \chi_{\{Q(s)=x\}}\mathcal{T}(ds|x',u')\nu_{i}(dx')\geq\int_{B_{i}}\int \chi_{\{Q(s)=x\}}\zeta(ds)\nu_{i}(dx')= \zeta(B_{i})>0   
\end{eqnarray}

Hence, every state action pair  $(X_{i},U_{i})\in \mathbb{X}^{m}\times \mathbb{U}^{m}$ is visited infinitely often. Thus, since under the minorization Assumption \ref{assumption for quantizing actions} the induced Markov chain is positive Harris recurrent, it follows from the ergodic theorem for positive Harris recurrent Markov chains (\cite[Theorem 17.1.7]{MeynBook} or \cite[Theorem 4.2.13]{HernandezLasserreErgodic}), that for all $\displaystyle (i,j,k)\in\{1,...,n\}^{3}$: \[M_{T}(x_{k}|x_{i},x_{j})=\frac{1}{\sum_{m=0}^{T-1}\chi_{\{X_{m}=x_{i},U_{m}=u_{j}\}}}\sum_{l=1}^{T}\chi_{\{X_{l}=x_{k},(X_{l-1}=x_{i},U_{l-1}=u_{j})\}}\rightarrow P_{n}(x_{k}|x_{i},u_{j})\] almost surely where $P_{n}$ is as in (\ref{P_{n}}) and $\nu$ is the invariant probability measure, which is induced by selecting actions from $\mathbb{U}^{m}$ uniformly at random, and satisfies $\nu(B_{i})>0$ for all bins $B_{i}$. Thus, almost surely, $M_{T}(\omega) \rightarrow P_{n}$.\\

\noindent Step 2: Let $v:\mathbb{X}\rightarrow \mathbb{R}$. We have that  
\begin{align*}
& \mathbb{T}^{T}v(x) :=\min_{u \in \mathbb{U}^{m}} \bigg( c(x,u)+\sum_{y\in \mathbb{X}^{m}}v(y)M_{T}(y|x,u)\bigg)\rightarrow \min_{u \in \mathbb{U}^{m}} \bigg( c(x,u)+\sum_{y\in \mathbb{X}^{m}}v(y)P_{n}(y|x,u)\bigg)
\end{align*} 

Following as in \cite{ky2026qaverage} via contraction of the average cost optimality operator obtained for the sub-stochastic kernel ${\cal T}(dX_1|x,u) - \zeta(dx)$, it follows that an optimal policy for the finite models exist and also for $T$ sufficiently large with $\gamma(M_{T}):=\argmin_{u \in \mathbb{U}^{m}}\bigg( c(x,u)+\sum_{y\in \mathbb{X}^{m}}v(y)M_{T}(y|x,u)\bigg)$, $\gamma(M_{T})\rightarrow \gamma(P_{n})$ is optimal for the finite state MDP whose transition kernel is given by $P_{n}$. 


\noindent Step 3: Then, we get $|J(\gamma_{T}^{*}) -J^{*}(\mathcal{T})|\leq |J(\gamma_{T}^{*}) -J^{*}_{n}| +|J^{*}_{n}-J^{*}(\mathcal{T})|$. It follows from Step 2 that the first term in the last inequality goes to zero as $T\rightarrow \infty$ whereas it follows from Theorem \ref{action quantization theorem} and \ref{theorem sate space quantization} that the second term goes to zero as $n\rightarrow \infty$.
\qed 


\begin{remark}\label{quantizedConvContConvremark}
Towards studying Algorithm \ref{alg:alg1}, it is important to note that the convergence result in Theorem \ref{theorem sate space quantization} is compatible with Theorem \ref{Theorem 2 Ali} in the sense that the extension of the quantized kernels $P_{n}$, which we denote $\Bar{P_{n}}$, by approximating inputs using the closest bin representative and expanding the image using uniform selection from each bin, leads to sequence of kernels that converges weakly continuously to ${\cal T}$ by \cite[Section 4]{KaraYuksel2021Chapter} and therefore the sequence converges uniformly by Theorem \ref{kernel topologies} (ii). As will be shown below $M_{T}(\omega) \rightarrow P_{n}$ pointwise and, due to finiteness, uniformly over all elements of $\mathbb{X}^{m}\times \mathbb{U}^{m}$. Accordingly, one obtains the following
 $\rho^{\mathrm{Unif}}_{BL}({\cal T},\Bar{M_T}(\omega))\leq \rho^{\mathrm{Unif}}_{BL}({\cal T},\Bar{P_n}) + \rho^{\mathrm{Unif}}_{BL}(\Bar{P_n},\Bar{M_T}(\omega))$, where $\Bar{M_T}(\omega)$ is the extension of $M_T(\omega)$. The first error term is a deterministic quantity and the second is a random quantity. Almost surely, for every $\omega$ and $\epsilon > 0$ $\exists n,T \in \mathbb{N}$ so that $\rho^{\mathrm{Unif}}_{BL}(\Bar{P_n},\Bar{M_T}(\omega)) \leq \epsilon$. This, then, leads to convergence to a near optimal solution via Theorem \ref{Theorem 2 Ali}.
\end{remark}

{
\begin{remark} [On Compactness]
    As can be seen from the results of Section 2 and 3, compactness of $\mathbb{X}$ was not necessary to achieve Bayesian identifiability and learning. In this section, we imposed compactness of $\mathbb{X}$ to avoid additional technical conditions and to preserve the same level of generality required on the kernels and cost functions. We note, however, that  approximation of Markov decision processes with $\sigma$-compact state space using an MDP with compact state space has been shown to be near optimal in \cite[Theorem 4.34]{SaLiYuSpringer} and \cite{SaYuLi15e} when one imposes additional conditions. Under these, one can first approximate $\mathbb{X}$ using a compact subset and then apply the results of this section. A further approach which leads to the quantization of MDPs with a $\sigma$-compact state space is presented in \cite{bicer2025quantizer} and shown to be near optimal under a Lyapunov condition and Lipschitz continuity of both the transition kernel and cost functional. 
\end{remark}}

\subsection{Numerical simulation}
\begin{figure}[h!]  
\centering
\includegraphics[width=0.7\columnwidth]{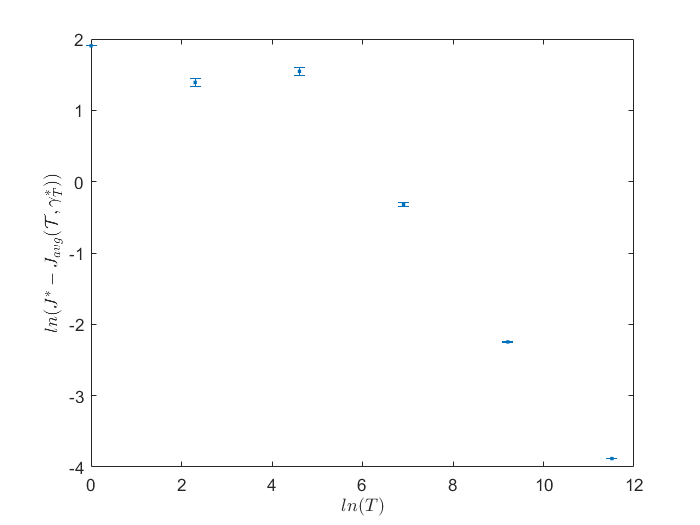}
\caption{$Err(T)$ vs $ln(T)$}
\label{fig:my_plot} 
\end{figure}

{We apply Algorithm I to the discrete time queuing problem provided in \cite[Section VI]{Kim2017TSfinite-case-average-cost}. Let $\gamma^{*}_{T}$ denote the optimal policy which is produced by running Algorithm I for horizon $T$. In what follows we evaluate $\frac{1}{10^{6}}\sum_{t=0}^{10^{6}}r(X_{t},U_{t})$, where $r(X_{t},U_{t})$ is the one stage reward function, when the policy $\gamma^{*}_{T}$ is applied and then compute the average performance across $1000$ samples. We denote the latter by $J_{avg}(\mathcal{T},\gamma^{*}_{T})$. Figure \ref{fig:my_plot} shows the difference between the optimal performance given the true kernel that is $J^{*}=1.9091$, and $J_{avg}(\mathcal{T},\gamma^{*}_{T})$ for various values of $T$. Here, $Err(T)=|J^{*}-J_{avg}(\mathcal{T},\gamma^{*}_{T})|$.}

\section{Concluding Remarks}

We considered a Markov Decision process with standard Borel spaces and an unknown transition kernel under an average cost criterion. We presented and studied several topologies on stochastic kernels which lead to their identifiability and robustness. We then obtained two data-driven identifiability results; the first one being Bayesian and the second one empirical. These were utilized to design near-optimal adaptive control policies with rigorous convergence guarantees.


We note that the Bayesian learning results required more stringent conditions. The Bayesian adaptive control techniques presented in this paper may be preferable in cases where only limited data can be gathered. The adaptive empirical learning, by contrast, does not require any restrictive regularity assumptions. Another advantage is that it does not require a prior on the set of kernels or the computation of Bayesian updates. 


\section*{Appendix}
\section*{Proof of Lemma \ref{M is tail measurable}}

Under the assumptions of Lemma \ref{M is tail measurable} we have that $\mathbb{M}$ is compact under the uniform bounded Lipschitz metric. Hence, it follows from \cite[Theorem 16.11]{Willard1970GeneralTopology} that $\mathbb{M}$ is separable. By \cite[Theorem 29]{DruryMATH354} this implies that every open set in $\mathbb{M}$ can be written as the countable union of open balls. Thus, to show that $M$ is $\mathcal{F}_{\infty}$ measurable, it is sufficient to show that for any $\epsilon>0$ and any $M_{j}\in \mathbb{M}$: $ \{M\in Q_{j} \}\in \mathcal{F}_{\infty}$ where $Q_{j}=\{m\in \mathbb{M}| \rho^{u}_{BL}(m,M_{j})\leq \epsilon\}$.

To this end let $\epsilon>0$. By the ergodic theorem for positive Harris recurrent Markov chains, we have that for all $B$ such that $\phi(B)>0$ and for any measurable set $A$ the following holds: 

\begin{align*}
   & P_{T-1}(A|B):=\frac{1}{\sum_{t=0}^{T-1}\chi_{\{(X_{t},U_{t})\in B\}}} \sum_{t=0}^{T-1}\chi_{\{X_{t+1}\in A , (X_{t},U_{t})\in B\}}\rightarrow \int_{B}\mathcal{T}(A|x,u)\Bar{\pi}(dx,du)
\end{align*}
$P^{\mathcal{T}}$ almost surely as $T\rightarrow \infty$.

Where, $P^{\mathcal{T}}$ is the measure induced on $\{(X_{t},U_{t})\}_{t=0}^{\infty}$ under the transition kernel $\mathcal{T}$, and $\Bar{\pi}$ is the unique invariant probability measure of the process $(X_{t},U_{t})$. Thus, for any pair $(B,A)$, $\displaystyle \lim_{T\rightarrow \infty} P_{T-1}(A|B)$ exists almost surely and thus, by definition of the tail $\sigma-$algebra, it follows that $\displaystyle \lim_{T\rightarrow \infty} P_{T-1}(A|B)$ is $\mathcal{F}_{\infty}$-measurable. First assume that $\mathbb{X}\times \mathbb{U}$ is compact. Recall that $\mathbf{F}$ is a countable collection of bounded Lipschitz functions which characterize the weak topology (Definition \ref{F definition}).

\begin{lemma}
We have that \begin{eqnarray*}
 \nonumber   \{M\in Q_{j}\}=\bigcap_{k=1}^{\infty}\bigcap_{f\in \mathbf{F}} \bigcap_{(\Bar{x}_{i},\Bar{u}_{i})\in S_{k}}\bigg\{\bigg|   \int f(s) M(ds|\Bar{x}_{i},\Bar{u}_{i})- \int f(s) M_{j}(ds|\Bar{x}_{i},\Bar{u}_{i})  \bigg| \leq \frac{1}{k} +\epsilon \bigg\}
\end{eqnarray*}
\noindent Where, $S_{k}$ is a quantization of $\mathbb{X}\times\mathbb{U}$ such that for all $(\Bar{x}_{i},\Bar{u}_{i})\in S_{k}$ $\exists$ $ \delta_{k}>0$ such that for all $(x,u)\in \mathbb{X}\times\mathbb{U}$: $d((x,u),(\Bar{x}_{i},\Bar{u}_{i}))<\delta_{k} \implies \forall M \in \mathbb{M}$: $\rho(M(.|x,u),M(.|\Bar{x}_{i},\Bar{u}_{i}))\leq \frac{1}{2(k+1)\|f\|_{BL}}$.
\end{lemma}

 We note that the existence of $S_{k}$ follows from Assumption \ref{equicontionuity}. Without loss of generality, one can choose the values of $\delta_{k}$ so that the sequence $\{\delta_{k}\}_{k\in \mathrm{N}}$ is decreasing, and for every $f \in \mathbf{F}$: $\delta_{k}\rightarrow 0$ as $k\rightarrow \infty$. In particular, $S_{k+1}$ forms a finer partition of $\mathbb{X}\times\mathbb{U}$ than $S_{k}$ in the sense that for all $k \in \mathbb{N}$, $S_{k}\subseteq S_{k+1}$.

\noindent Proof.
Suppose for some $f\in \mathbf{F}$ that the following holds

\noindent $\displaystyle \bigcap_{k=1}^{\infty}\bigcap_{(\Bar{x}_{i},\Bar{u}_{i})\in S_{k}}\bigg\{\bigg|  \int f(s) M(ds|\Bar{x}_{i},\Bar{u}_{i})- \int f(s) M_{j}(ds|\Bar{x}_{i},\Bar{u}_{i})  \bigg|\leq \frac{1}{k}+\epsilon \bigg\}$. Then, we get that
\begin{eqnarray*}
\nonumber    &&\bigcap_{(x,u)\in \mathbb{X}\times\mathbb{U}} \bigg\{\bigg| \int f(s) M(ds|x,u)- \int f(s) M_{j}(ds|x,u)  \bigg| \leq \epsilon \bigg\}=\\
\nonumber    &&\bigcap_{k=1}^{\infty}\bigcap_{(x,u)\in \mathbb{X}\times\mathbb{U}} \bigg\{\bigg|   \int f(s) M(ds|x,u)- \int f(s) M_{j}(ds|x,u)  \bigg| \leq \frac{1}{k}+\epsilon \bigg\}
\end{eqnarray*}
Suppose that for all $k\in \mathbb{N}$ the following holds for all $\Bar{x}_{i},\Bar{u}_{i}\in S_{k}$ 

 \[\bigg| \int f(s) M(ds|\Bar{x}_{i},\Bar{u}_{i})- \int f(s) M_{j}(ds|\Bar{x}_{i},\Bar{u}_{i})  \bigg| \leq \frac{1}{k}+\epsilon\]
 
 Then,  for any $(x,u)\in \mathbb{X}\times\mathbb{U}$
\begin{eqnarray*}
\nonumber    \bigg|   \int f(s) M(ds|x,u)- \int f(s) M_{j}(ds|x,u)  \bigg| &\leq& \bigg|   \int f(s) M(ds|x,u)- \int f(s) M(ds|\Bar{x}_{i},\Bar{u}_{i})  \bigg|\\
 \nonumber   && + \quad \bigg|   \int f(s) M(ds|\Bar{x}_{i},\Bar{u}_{i})- \int f(s) M_{j}(ds|\Bar{x}_{i},\Bar{u}_{i})  \bigg|\\
\nonumber    &&+ \quad \bigg|   \int f(s) M_{j}(ds|\Bar{x}_{i},\Bar{u}_{i})- \int f(s) M_{j}(ds|x,u)  \bigg|\\
&&\leq \frac{1}{2(k+1)}+ \frac{1}{k} + \epsilon+ \frac{1}{2(k+1)}=\frac{1}{k} + \epsilon+ \frac{1}{k+1}
\end{eqnarray*}

Since the last inequality holds for all $k\in \mathbb{N}$, it follows that for any $k'\in \mathbb{N}$ 
\begin{eqnarray*}
\nonumber        \bigg|  \int f(s) M(ds|x,u)- \int f(s) M_{j}(ds|x,u)  \bigg| &\leq& \frac{1}{k'}+\epsilon
\end{eqnarray*}
Thus we have shown that 
\begin{eqnarray*}
 \nonumber    \{M\in Q_{j}\}&=& \bigcap_{f\in \mathbf{F}}\bigcap_{k=1}^{\infty}\bigcap_{(x,u)\in \mathbb{X}\times\mathbb{U}} \bigg\{\bigg|  \int f(s) M(ds|x,u)- \int f(s) M_{j}(ds|x,u)  \bigg| \leq \frac{1}{k}+\epsilon \bigg\}\\
\nonumber     &\supseteq&\bigcap_{f\in \mathbf{F}}\bigcap_{k=1}^{\infty} \bigcap_{(\Bar{x}_{i},\Bar{u}_{i})\in S_{k}}\bigg\{\bigg|   \int f(s) M(ds|\Bar{x}_{i},\Bar{u}_{i})- \int f(s) M_{j}(ds|\Bar{x}_{i},\Bar{u}_{i})  \bigg| \leq \frac{1}{k} +\epsilon \bigg\}
\end{eqnarray*}
Since, the other inclusion holds trivially, this completes the proof of the claim. The following lemma will also be essential to complete the proof.

\begin{lemma}
\begin{eqnarray*}
 \nonumber   &&\bigcap_{k=1}^{\infty}\bigcap_{f\in \mathbf{F}} \bigcap_{(\Bar{x}_{i},\Bar{u}_{i})\in S_{k}}\bigg\{\bigg|   \int f(s) M(ds|\Bar{x}_{i},\Bar{u}_{i})- \int f(s) M_{j}(ds|\Bar{x}_{i},\Bar{u}_{i})  \bigg| \leq \frac{1}{k}+\epsilon  \bigg\}=\\
 \nonumber   &&\bigcap_{k=1}^{\infty}\bigcap_{f\in \mathbf{F}} \bigcap_{(\Bar{x}_{i},\Bar{u}_{i})\in S_{k}}\bigg\{\bigg|   \int_{B_{\delta}(\Bar{x}_{i},\Bar{u}_{i})}\int f(s) M(ds|x,u)\Bar{\pi}(dx,du)\\
 \nonumber &&- \quad \int_{B_{\delta}(\Bar{x}_{i},\Bar{u}_{i})}\int f(s) M_{j}(ds|x,u)\Bar{\pi}(dx,du) \bigg| \leq \frac{1}{k} +\epsilon \bigg\}
\end{eqnarray*}  
\end{lemma}
Proof.
Let $f \in \mathbf{F}$. Suppose that for all $k\in \mathrm{N}$, we have that \\
\[\bigg|   \int f(s) M(ds|\Bar{x}_{i},\Bar{u}_{i})- \int f(s) M_{j}(ds|\Bar{x}_{i},\Bar{u}_{i})  \bigg| \leq \frac{1}{k}+\epsilon \]

Let $K\in \mathbb{N}$. Then,

\begin{eqnarray*}
\nonumber   && \bigg|   \int_{B_{\delta}(\Bar{x}_{i},\Bar{u}_{i})}\int f(s) M(ds|x,u)\Bar{\pi}(dx,du)- \int_{B_{\delta}(\Bar{x}_{i},\Bar{u}_{i})}\int f(s) M_{j}(ds|x,u)\Bar{\pi}(dx,du) \bigg|\\
\nonumber   &&\leq  \bigg|  \int_{B_{\delta}(\Bar{x}_{i},\Bar{u}_{i})}\int f(s) M(ds|x,u)\Bar{\pi}(dx,du)- \int_{B_{\delta}(\Bar{x}_{i},\Bar{u}_{i})}\int f(s) M(ds|\Bar{x}_{i},\Bar{u}_{i})\Bar{\pi}(dx,du) \bigg|  \\
\nonumber   &&+ \quad \bigg | \int_{B_{\delta}(\Bar{x}_{i},\Bar{u}_{i})}\int f(s) M(ds|\Bar{x}_{i},\Bar{u}_{i})\Bar{\pi}(dx,du)-\int_{B_{\delta}(\Bar{x}_{i},\Bar{u}_{i})}\int f(s) M_{j}(ds|\Bar{x}_{i},\Bar{u}_{i})\Bar{\pi}(dx,du) \bigg | +\\
\nonumber   &&+ \quad \bigg | \int_{B_{\delta}(\Bar{x}_{i},\Bar{u}_{i})}\int f(s) M_{j}(ds|\Bar{x}_{i},\Bar{u}_{i})\Bar{\pi}(dx,du)- \int_{B_{\delta}(\Bar{x}_{i},\Bar{u}_{i})}\int f(s) M_{j}(ds|x,u)\Bar{\pi}(dx,du) \bigg | \\
 \nonumber  &&\leq \frac{1}{k}+\frac{1}{k+1}+\epsilon \leq \frac{1}{K}+\epsilon 
\end{eqnarray*}
for some $k$ sufficiently large, which we denote by $k(K)$, so that the last inequality holds. Thus, because $S_{K}\subseteq S_{k}$, we get that 
\begin{align} \label{inclusion uno}
    &\bigcap_{(\Bar{x}_{i},\Bar{u}_{i})\in S_{k}}\bigg\{\bigg |   \int f(s) M(ds|\Bar{x}_{i},\Bar{u}_{i})- \int f(s) M_{j}(ds|\Bar{x}_{i},\Bar{u}_{i})  \bigg | \leq \frac{1}{k} +\epsilon \bigg\} \subseteq \\
   \nonumber  & \bigcap_{(\Bar{x}_{i},\Bar{u}_{i})\in S_{K}}\bigg\{\bigg |   \int_{B_{\delta}(\Bar{x}_{i},\Bar{u}_{i})}\int f(s) M(ds|x,u)\Bar{\pi}(dx,du)-\int_{B_{\delta}(\Bar{x}_{i},\Bar{u}_{i})}\int f(s) M_{j}(ds|x,u)\Bar{\pi}(dx,du) \bigg | \leq \frac{1}{K}+\epsilon  \bigg\} 
\end{align}
Since the inclusion (\ref{inclusion uno}) holds for every $f\in \mathbf{F}$, it follows that 
\begin{eqnarray*}
\nonumber    &&\bigcap_{f\in \mathbf{F}}\bigcap_{(\Bar{x}_{i},\Bar{u}_{i})\in S_{k}}\bigg\{\bigg |   \int f(s) M(ds|\Bar{x}_{i},\Bar{u}_{i})- \int f(s) M_{j}(ds|\Bar{x}_{i},\Bar{u}_{i})  \bigg | \leq \frac{1}{k} +\epsilon \bigg\} \subseteq \\
\nonumber    && \bigcap_{f\in \mathbf{F}}\bigcap_{(\Bar{x}_{i},\Bar{u}_{i})\in S_{K}}\bigg\{\bigg |   \int_{B_{\delta}(\Bar{x}_{i},\Bar{u}_{i})}\int f(s) M(ds|x,u)\Bar{\pi}(dx,du) \\
\nonumber &&- \quad \int_{B_{\delta}(\Bar{x}_{i},\Bar{u}_{i})}\int f(s) M_{j}(ds|x,u)\Bar{\pi}(dx,du) \bigg | \leq \frac{1}{K} +\epsilon \bigg\}
\end{eqnarray*}
Thus, we get that 
\begin{eqnarray*}
 \nonumber   &&\bigcap_{K=1}^{\infty}\bigcap_{f\in \mathbf{F}}\bigcap_{(\Bar{x}_{i},\Bar{u}_{i})\in S_{K}}\bigg\{\bigg |   \int f(s) M(ds|\Bar{x}_{i},\Bar{u}_{i})- \int f(s) M_{j}(ds|\Bar{x}_{i},\Bar{u}_{i})  \bigg | \leq \frac{1}{K}+\epsilon  \bigg\} \subseteq \\
\nonumber    &&\bigcap_{K=1}^{\infty}\bigcap_{k(K)}\bigcap_{f\in \mathbf{F}}\bigcap_{(\Bar{x}_{i},\Bar{u}_{i})\in S_{k}}\bigg\{\bigg |   \int f(s) M(ds|\Bar{x}_{i},\Bar{u}_{i})- \int f(s) M_{j}(ds|\Bar{x}_{i},\Bar{u}_{i})  \bigg | \leq \frac{1}{k} +\epsilon \bigg\} \subseteq\\
\nonumber    && \bigcap_{K=1}^{\infty}\bigcap_{f\in \mathbf{F}}\bigcap_{(\Bar{x}_{i},\Bar{u}_{i})\in S_{K}}\bigg\{\bigg |   \int_{B_{\delta}(\Bar{x}_{i},\Bar{u}_{i})}\int f(s) M(ds|x,u)\Bar{\pi}(dx,du)\\
\nonumber &&- \quad \int_{B_{\delta}(\Bar{x}_{i},\Bar{u}_{i})}\int f(s) M_{j}(ds|x,u)\Bar{\pi}(dx,du) \bigg | \leq \frac{1}{K} +\epsilon \bigg\}
\end{eqnarray*}

Now, suppose for some $k\in \mathrm{N}$, $f \in \mathbf{F}$, and $(\Bar{x}_{i},\Bar{u}_{i})\in S_{k}$ that 
\[\bigg |   \int f(s) M(ds|\Bar{x}_{i},\Bar{u}_{i})- \int f(s) M_{j}(ds|\Bar{x}_{i},\Bar{u}_{i}) \bigg | > \frac{1}{k} +\epsilon \] 

Without loss of generality, we can assume that \[\int f(s) M(ds|\Bar{x}_{i},\Bar{u}_{i})-\int f(s) M_{j}(ds|\Bar{x}_{i},\Bar{u}_{i}) > \frac{1}{k}+\epsilon \]

By the continuity of $\int f(s) M(ds|x,u)- \int f(s) M_{j}(ds|x,u)$ in $(x,u)$, it follows that for some $K$ large enough we have that for all $(x,u) \in B_{\delta_{K}}$: $\int f(s) M(ds|x,u)- \int f(s) M_{j}(ds|x,u)>\frac{1}{K} +\epsilon $ and thus

\begin{eqnarray*}
\nonumber    &&\Big |   \int_{B_{\delta}(\Bar{x}_{i},\Bar{u}_{i})}\int f(s) M(ds|x,u)\Bar{\pi}(dx,du)- \int_{B_{\delta}(\Bar{x}_{i},\Bar{u}_{i})}\int f(s) M_{j}(ds|x,u)\Bar{\pi}(dx,du) \Big | =\\
\nonumber    && \Big |   \int_{B_{\delta}(\Bar{x}_{i},\Bar{u}_{i})}\Bar{\pi}(dx,du)\Big[\int f(s) M(ds|x,u)- \int f(s) M_{j}(ds|x,u) \Big] \Big |>\frac{1}{K}+\epsilon 
\end{eqnarray*}

Hence,
\begin{eqnarray*}
\nonumber    &&\Bigg (\bigcap_{k=1}^{\infty}\bigcap_{f\in \mathbf{F}} \bigcap_{(\Bar{x}_{i},\Bar{u}_{i})\in S_{k}}\bigg\{\bigg |   \int f(s) M(ds|\Bar{x}_{i},\Bar{u}_{i})- \int f(s) M_{j}(ds|\Bar{x}_{i},\Bar{u}_{i})  \bigg | \leq \frac{1}{k} +\epsilon \bigg\}\Bigg )^{C} \subseteq\\
\nonumber    &&\Bigg (\bigcap_{k=1}^{\infty}\bigcap_{f\in \mathbf{F}} \bigcap_{(\Bar{x}_{i},\Bar{u}_{i})\in S_{k}}\bigg\{\bigg |   \int_{B_{\delta}(\Bar{x}_{i},\Bar{u}_{i})}\int f(s) M(ds|x,u)\Bar{\pi}(dx,du)- \\
\nonumber &&\int_{B_{\delta}(\Bar{x}_{i},\Bar{u}_{i})}\int f(s) M_{j}(ds|x,u)\Bar{\pi}(dx,du) \bigg | \leq \frac{1}{k}+\epsilon  \bigg\}\Bigg )^{C}
\end{eqnarray*}

which implies that

\begin{eqnarray*}
 \nonumber   &&\bigcap_{k=1}^{\infty}\bigcap_{f\in \mathbf{F}} \bigcap_{(\Bar{x}_{i},\Bar{u}_{i})\in S_{k}}\bigg\{\bigg |   \int_{B_{\delta}(\Bar{x}_{i},\Bar{u}_{i})}\int f(s) M(ds|x,u)\Bar{\pi}(dx,du) - \quad \int_{B_{\delta}(\Bar{x}_{i},\Bar{u}_{i})}\int f(s) M_{j}(ds|x,u)\Bar{\pi}(dx,du) \bigg | \leq \frac{1}{k}+\epsilon  \bigg\}\\
\nonumber    &&\subseteq\bigcap_{k=1}^{\infty}\bigcap_{f\in \mathbf{F}} \bigcap_{(\Bar{x}_{i},\Bar{u}_{i})\in S_{k}}\bigg\{\bigg |   \int f(s) M(ds|\Bar{x}_{i},\Bar{u}_{i})- \int f(s) M_{j}(ds|\Bar{x}_{i},\Bar{u}_{i})  \bigg | \leq \frac{1}{k}+\epsilon  \bigg\}
\end{eqnarray*}
and this completes the proof of the Lemma. 

Hence, to show that $\{M\in Q_{j}\}\in \mathcal{F}_{\infty}$ it is sufficient to show that 

\begin{eqnarray} \label{conditionn}
   \nonumber &&\bigcap_{k=1}^{\infty}\bigcap_{f\in \mathbf{F}}\bigcap_{(\Bar{x}_{i},\Bar{u}_{i})\in S_{k}}\bigg\{\bigg |   \int_{B_{\delta}(\Bar{x}_{i},\Bar{u}_{i})}\int f(s) M(ds|x,u)\Bar{\pi}(dx,du)\\
    &&- \quad \int_{B_{\delta}(\Bar{x}_{i},\Bar{u}_{i})}\int f(s) M_{j}(ds|x,u)\Bar{\pi}(dx,du) \bigg | \leq \frac{1}{k}+\epsilon \bigg\}\in \mathcal{F}_{\infty}
\end{eqnarray}
Let $f\in \mathbf{F}$. Since $f$ is bounded, there exists a sequence of functions $S_{n}(w)=\sum_{l=0}^{n}\alpha_{l}^{n}\chi_{A^{n}_{l}}(w)$ which is monotonically increasing, uniformly bounded from below, and satisfies $S_{n}\rightarrow f$ pointwise. Thus, one can write for any $f$ and $(\Bar{x}_i,\Bar{u}_i)$

\begin{eqnarray*}
\nonumber    &&\int_{B_{\delta}(\Bar{x}_{i},\Bar{u}_{i})}\int f(s) M(ds|x,u)\Bar{\pi}(dx,du)=\int_{B_{\delta}(\Bar{x}_{i},\Bar{u}_{i})}\lim_{n\rightarrow \infty}\sum_{l=0}^{n}\alpha_{l}^{n}\mathcal{T}(A^{n}_{l}|x,u)\Bar{\pi}(dx,du)\\
\nonumber    &&=\lim_{n\rightarrow \infty}\int_{B_{\delta}(\Bar{x}_{i},\Bar{u}_{i})}\sum_{l=0}^{n}\alpha_{l}^{n}\mathcal{T}(A^{n}_{l}|x,u)\Bar{\pi}(dx,du)=\lim_{n\rightarrow \infty}\sum_{l=0}^{n}\alpha_{l}^{n}\int_{B_{\delta}(\Bar{x}_{i},\Bar{u}_{i})}\mathcal{T}(A^{n}_{l}|x,u)\Bar{\pi}(dx,du)\\
\nonumber    &&=\lim_{n\rightarrow \infty}\sum_{l=0}^{n}\alpha_{l}^{n}\lim_{T\rightarrow \infty}P_{T-1}(A^{n}_{l}|B_{\delta}
    (\Bar{x}_{i},\Bar{u}_{i}))
\end{eqnarray*}
 Because for any pair $(B,A)$, $\displaystyle \lim_{T\rightarrow \infty} P_{T-1}(A|B)$ is $\mathcal{F}_{\infty}$-measurable, it then follows that condition (\ref{conditionn}) holds and thus $\{M\in Q_{j}\}\in \mathcal{F}_{\infty}$. \qed

Now suppose $\mathbb{X}\times\mathbb{U}$ is not compact, but $\sigma$-compact. Then, one can write $\mathbb{X}\times\mathbb{U}=\bigcup_{n=1}^{\infty}E_{n}$ where for all $n$ $E_{n}$ is compact and $E_{n}\subseteq E_{n+1}$. Hence, we get that
\begin{eqnarray*}
\nonumber    \{M\in Q_{j}\}=\bigcap_{n=1}^{\infty}\bigcap_{k=1}^{\infty}\bigcap_{f\in \mathbf{F}} \bigcap_{(\Bar{x}_{i},\Bar{u}_{i})\in S_{k}^{n}}\bigg\{\bigg |   \int f(s) M(ds|\Bar{x}_{i},\Bar{u}_{i})- \int f(s) M_{j}(ds|\Bar{x}_{i},\Bar{u}_{i})  \bigg | \leq \frac{1}{k} +\epsilon\bigg\}
\end{eqnarray*}

\noindent Where, $S_{k}^{n}$ is a quantization of $E_{n}$ such that for all $(\Bar{x}_{i},\Bar{u}_{i})\in S_{k}$ $\exists$ $ \delta_{k}>0$ such that for all $(x,u)\in E_{n}$: $d((x,u),(\Bar{x}_{i},\Bar{u}_{i}))<\delta_{k} \implies \forall M \in \mathbb{M}$: $\rho(M(.|x,u),M(.|\Bar{x}_{i},\Bar{u}_{i}))\leq \frac{1}{2(k+1)\|f\|_{BL}}$. Because, for all $n$ we have that 
\begin{eqnarray*}
\nonumber    \bigcap_{k=1}^{\infty}\bigcap_{f\in \mathbf{F}} \bigcap_{(\Bar{x}_{i},\Bar{u}_{i})\in S_{k}^{n}}\bigg\{\bigg |   \int f(s) M(ds|\Bar{x}_{i},\Bar{u}_{i})- \int f(s) M_{j}(ds|\Bar{x}_{i},\Bar{u}_{i})  \bigg | \leq \frac{1}{k}+\epsilon \bigg\}\in \mathcal{F}_{\infty}
\end{eqnarray*}
it follows that $\{M\in Q_{j}\}\in \mathcal{F}_{\infty}\subseteq I_{\infty}:= \sigma(X_{0},U_{0},X_{1},U_{1},...)$ $P^{\mathcal{T}}$ almost surely where $P^{\mathcal{T}}$ denotes the measure induced on the process $\{(X_{t},U_{t})\}_{t=0}^{\infty}$ 

under the transition kernel $\mathcal{T}$. Since under any transition kernel $m\in \mathbb{M}$ the process $(X_{t},U_{t})$ is positive Harris recurrent and admits an exhaustive unique invariant probability measure, we get that for all $m\in\mathbb{M}$: $\{M\in {Q}_{j}\}\in \mathcal{F}_{\infty}\subseteq I_{\infty}:= \sigma(X_{0},U_{0},X_{1},U_{1},...)$ $P^{m}$ almost surely where $P^{m}$ denotes the measure induced on the process $(X_{t},U_{t})$ when the true kernel is $m\in \mathbb{M}$ and the exploration policy $\gamma$ is applied. Hence, $\{M\in {Q}_{j}\}\in I_{\infty}$ $P^{\Theta}$ almost surely. 
 \qed



\begin{thebibliography}{10}

\bibitem{abounadi2001learning}
J.~Abounadi, D.~Bertsekas, and V.S. Borkar.
\newblock Learning algorithms for {M}arkov decision processes with average cost.
\newblock {\em SIAM Journal on Control and Optimization}, 40(3):681--698, 2001.

\bibitem{aldous1981weak}
D.J. Aldous.
\newblock Weak convergence and the general theory of processes.
\newblock Incomplete draft of a monograph, July, 1981.

\bibitem{survey}
A.~Arapostathis, V.~S. Borkar, E.~Fernandez-Gaucherand, M.~K. Ghosh, and S.~I. Marcus.
\newblock Discrete-time controlled {M}arkov processes with average cost criterion: A survey.
\newblock {\em SIAM J. Control and Optimization}, 31:282--344, 1993.

\bibitem{julio2020adapted}
J.~Backhoff-Veraguas, D.~Bartl, M.~Beiglb{\"o}ck, and M.~Eder.
\newblock Adapted {W}asserstein distances and stability in mathematical finance.
\newblock {\em Finance and Stochastics}, 24(3):601--632, 2020.

\bibitem{backhoff2019all}
J.~Backhoff-Veraguas, D.~Bartl, M.~Beiglb{\"o}ck, and M.~Eder.
\newblock All adapted topologies are equal.
\newblock {\em Probability Theory and Related Fields}, 178(3):1125--1172, 2020.

\bibitem{BalderYoungMeasureTheory}
E.~J. Balder.
\newblock Lectures on young measure theory and its applications in economics.
\newblock {\em Rend. Ist. Mat. Univ. Trieste}, 31:1--69, 2000.

\bibitem{balder2001}
E.~J. Balder.
\newblock On ws-convergence of product measures.
\newblock {\em Mathematics of Operations Research}, 26(3):494--518, 2001.

\bibitem{Kim2019TS-continuous-case-average-cost}
D.~Banjevi{\'c} and M.J. Kim.
\newblock Thompson samplig for stochastic control: the continuous parameter case.
\newblock {\em IEE Transactions on Automatic Control}, 64:4137--4152, 2019.

\bibitem{bartl2024wasserstein}
D.~Bartl, M.~Beiglb{\"o}ck, and G.~Pammer.
\newblock The {W}asserstein space of stochastic processes.
\newblock {\em Journal of the European Mathematical Society}, 2024.

\bibitem{bayraktar2020continuity}
E.~Bayraktar, Y.~Dolinsky, and J.~Guo.
\newblock Continuity of utility maximization under weak convergence.
\newblock {\em Mathematics and Financial Economics}, pages 1--33, 2020.

\bibitem{beiglbock2022approximation}
M.~Beiglb{\"o}ck, B.~Jourdain, W.~Margheriti, and G.~Pammer.
\newblock Approximation of martingale couplings on the line in the adapted weak topology.
\newblock {\em Probability Theory and Related Fields}, 183(1):359--413, 2022.

\bibitem{bicer2025quantizer}
O.~Bi\c{c}er, A.D. Kara, and Serdar S.~Y\"uksel.
\newblock Quantizer design for finite model approximations, model learning, and quantized q-learning for mdps with unbounded spaces.
\newblock {\em arXiv preprint arXiv:2510.04355}, 2025.

\bibitem{Bor91}
V.~S. Borkar.
\newblock On extremal solutions to stochastic control problems.
\newblock {\em Applied Mathematics {\&} Optimization}, 24(1):317--330, 1991.

\bibitem{BorkarRealization}
V.~S. Borkar.
\newblock White-noise representations in stochastic realization theory.
\newblock {\em SIAM J. on Control and Optimization}, 31:1093--1102, 1993.

\bibitem{bozkurt2024modelapproximationmdpsunbounded}
B.~Bozkurt, A.~Mahajan, A.~Nayyar, and Y.~Ouyang.
\newblock Model approximation in mdps with unbounded per-step cost.
\newblock {\em IEEE Transactions on Automatic Control}, pages 1--16, 2025.

\bibitem{Ionescu-Tulcea}
C.Ionescu-Tulcea.
\newblock Mesures dans les espaces produits.
\newblock {\em Atti Accad. Naz. Lincei Rend}, 7:208--211, 1949.

\bibitem{CostaDufourGenadot2026}
O.~L.~V. Costa, F.~Dufour, and A.~Genadot.
\newblock Adaptive nonstationary value iteration for discounted control of piecewise deterministic markov processes.
\newblock {\em Mathematics of Control, Signals, and Systems}, 2026.
\newblock Published online: 30 April 2026.

\bibitem{Adelman2024Thompsonsampling}
C.Keceli D.Adelman and A.V.~Olivares Nadal.
\newblock Thompson sampling for inifinite-horizon discounted decision process.
\newblock {\em arXiv preprint arXiv:2405.08253}, 2024.

\bibitem{bartl2023sensitivity}
D.Bartl and J.Wiesel.
\newblock Sensitivity of multiperiod optimization problems with respect to the adapted wasserstein distance.
\newblock {\em SIAM Journal on Financial Mathematics}, 14(2):704--720, 2023.

\bibitem{Egorovtheorem1911}
D.F.Egorov.
\newblock Sur les suites des fonctions mesurables.
\newblock {\em Comptes rendus hebdomadaires des s\'eances de l'Acad\'emie des sciences}, 152:135--157, 1911.

\bibitem{DruryMATH354}
S.~W. Drury.
\newblock Honours analysis 3: {MATH} 354.
\newblock Lecture notes, McGill University.
\newblock Metric spaces and calculus in several variables.

\bibitem{DuPr14}
F.~Dufour and T.~Prieto-Rumeau.
\newblock Approximation of average cost {M}arkov decision processes using empirical distributions and concentration inequalities.
\newblock {\em Stochastics An International Journal of Probability and Stochastic Processes}, 87:273--307, 2015.

\bibitem{Gordienko1985}
E.I.Gordienko.
\newblock Adaptive strategies for certain classes of controlled markov processes.
\newblock {\em Theory of Probability and Its Applications}, 29:504--518, 1985.

\bibitem{ethier2009markov}
S.~N. Ethier and T.~G. Kurtz.
\newblock {\em Markov Processes: Characterization and Convergence}, volume 282.
\newblock John Wiley \& Sons, 2009.

\bibitem{On-the-certainty-equivalence}
F.D\"orfler, P.Tesi, and C.~De{ }Persis.
\newblock On the certainty-equivalence approach to direct data-driven {LQR} design.
\newblock {\em IEEE Transactions on Automatic Control}, 68:7989--7996, 2023.

\bibitem{Caines-adaptive-control}
G.C.Goodwin, P.J.Ramadge, and P.E.Caines.
\newblock Discrete time stochastic adaptive control.
\newblock {\em SIAM Journal On Control and Optimization}, 19:829--853, 1981.

\bibitem{gihman2012controlled}
I.~I. Gihman and A.~V. Skorohod.
\newblock {\em Controlled Stochastic Processes}.
\newblock Springer Science \& Business Media, 2012.

\bibitem{GordienkoODE}
E.~I. Gordienko.
\newblock An estimate of the stability of optimal control of certain stochastic and deterministic systems.
\newblock {\em Journal of Soviet Mathematics}, 59(4):891–899, April 1992.

\bibitem{GordienkoSystemControlLetters}
E.~I. Gordienko and F.~Salem-Silva.
\newblock Robustness inequality for {M}arkov control processes with unbounded costs.
\newblock {\em Systems and Control Letters}, 33(2):125–130, February 1998.

\bibitem{GordienkoKybernetic}
E.~I. Gordienko and F.~Salem-Silva.
\newblock Estimates of stability of {M}arkov control processes with unbounded costs.
\newblock {\em Kybernetika}, 36(2):[195]--210, 2000.

\bibitem{gosavi2004reinforcement}
A.~Gosavi.
\newblock Reinforcement learning for long-run average cost.
\newblock {\em European journal of operational research}, 155(3):654--674, 2004.

\bibitem{hellwig1996sequential}
M.~F. Hellwig.
\newblock Sequential decisions under uncertainty and the maximum theorem.
\newblock {\em Journal of Mathematical Economics}, 25(4):443--464, 1996.

\bibitem{hernandez2012adaptive}
O.~Hern{\'a}ndez-Lerma.
\newblock {\em Adaptive Markov control processes}, volume~79.
\newblock Springer Science \& Business Media, 2012.

\bibitem{HernandezLasserreErgodic}
O.~Hern{\'a}ndez-Lerma and J.~B. Lasserre.
\newblock {\em {M}arkov Chains and Invariant Probabilities}, volume 211.
\newblock Birkh\"auser, Basel, 2003.

\bibitem{bayrooti2025noregret}
J.Bayrooti, S.Vakili, A.Prorok, and C.Henrik EK.
\newblock No-regret {T}hompson sampling for finite-horizon markov decision processes with gaussian processes.
\newblock In {\em Proceedings of the 39th Conference on Neural Information Processing Systems}, 2025.

\bibitem{kara2022robustness}
A.~D. Kara, M.~Raginsky, and S.~Y{\"u}ksel.
\newblock Robustness to incorrect models and data-driven learning in average-cost optimal stochastic control.
\newblock {\em Automatica}, 139:110179, 2022.

\bibitem{kara2020robustness}
A.D Kara and S.~Y\"uksel.
\newblock Robustness to incorrect system models in stochastic control.
\newblock {\em SIAM Journal on Control and Optimization}, 58(2):1144--1182, 2020.

\bibitem{KaraYuksel2021Chapter}
A.D. Kara and S.~Y\"uksel.
\newblock Robustness to approximations and model learning in mdps and pomdps.
\newblock In A.~B. Piunovskiy and Y.~Zhang, editors, {\em Modern Trends in Controlled Stochastic Processes: Theory and Applications, Volume III}. Luniver Press, 2021.

\bibitem{ky2026qaverage}
A.D. Kara and S.~Y\"uksel.
\newblock Approximations and learning for continuous state and action mdps under average cost criteria.
\newblock {\em Journal of Machine Learning Research}, 2026.

\bibitem{Kim2017TSfinite-case-average-cost}
M.J. Kim.
\newblock Thompson samplig for stochastic control: the finite parameter case.
\newblock {\em IEEE Transactions on Automatic Control}, 62:6415--6422, 2017.

\bibitem{Lan81}
H.J. Langen.
\newblock Convergence of dynamic programming models.
\newblock {\em Mathematics of Operations Research}, 6(4):493--512, Nov. 1981.

\bibitem{MeynBook}
S.~P. Meyn and R.~Tweedie.
\newblock {\em {M}arkov Chains and Stochastic Stability}.
\newblock Springer-Verlag, London, 1993.

\bibitem{muller1997does}
A.~M{\"u}ller.
\newblock How does the value function of a markov decision process depend on the transition probabilities?
\newblock {\em Mathematics of Operations Research}, 22(4):872--885, 1997.

\bibitem{munkres2000topology}
J.~R. Munkres.
\newblock {\em Topology}.
\newblock Prentice Hall, Upper Saddle River, NJ, 2 edition, 2000.

\bibitem{Zurek-plugin-approach}
M.Zurek and Y.Chen.
\newblock The plug-in approach for average-reward and discounted mdps: Optimal sample complexity analysis.
\newblock {\em arXiv preprint arXiv:2410.07616}, 2024.

\bibitem{Hernandez-lerma1987AdaptiveControl}
O.Hern\'andez-Lerma.
\newblock Approximation and adaptive control of markov processes: Average reward criterion.
\newblock {\em Kybernetika}, 23:265--288, 1987.

\bibitem{fatou-lemma}
O.Hernandez-Lerma and J.B.Lasserre.
\newblock Fatou’s lemma and {L}ebesgue’s convergence theorem for measures.
\newblock {\em Journal of Applied Mathematics and Stochastic Analysis}, 13(2):137--146, 2000.

\bibitem{Hernandez-Lerma1990DensityestimationAndAdaptiveControl}
O.Hern\'andez-Lerma and R.Cavazos-Cadena.
\newblock Density estimation and adaptive control of markov processes: Average and discounted criteria.
\newblock {\em Acta Applicandae Mathematica}, 20:285--307, 1990.

\bibitem{pammer2024note}
G.~Pammer.
\newblock A note on the adapted weak topology in discrete time.
\newblock {\em Electronic Communications in Probability}, 29:1--13, 2024.

\bibitem{pradhan2022robustness}
S.~Pradhan and S.~Y\"uksel.
\newblock Robustness of stochastic optimal control to approximate diffusion models under several cost evaluation criteria.
\newblock {\em Mathematics of Operations Research}, 2023.

\bibitem{russo2018tutorial}
D.J. Russo, B.~Van Roy, A.~Kazerouni, I.~Osband, and Z.~Wen.
\newblock A tutorial on thompson sampling.
\newblock {\em Foundations and Trends in Machine Learning}, 11(1):1--96, 2018.

\bibitem{SaLiYuSpringer}
N.~Saldi, T.~Linder, and S.~Y\"uksel.
\newblock {\em Finite Approximations in Discrete-Time Stochastic Control: Quantized Models and Asymptotic Optimality}.
\newblock Springer, Cham, 2018.

\bibitem{saldi2025kernel}
N.~Saldi and S.~Y\"uksel.
\newblock Kernel mean embedding topology: Weak and strong forms for stochastic kernels and implications for model learning.
\newblock {\em Annals of Applied Probability}, 2026.

\bibitem{SaYuLi15e}
N.~Saldi, S.~Y\"uksel, and T.~Linder.
\newblock Finite-state approximation of {M}arkov decision processes with unbounded costs and {B}orel spaces.
\newblock In {\em IEEE Conf. Decision Control}, Osaka. Japan, December 2015.

\bibitem{SaYuLi15c}
N.~Saldi, S.~Y{\"u}ksel, and T.~Linder.
\newblock On the asymptotic optimality of finite approximations to {M}arkov decision processes with {B}orel spaces.
\newblock {\em Mathematics of Operations Research}, 42(4):945--978, 2017.

\bibitem{Schal}
M.~Sch\"al.
\newblock Conditions for optimality in dynamic programming and for the limit of n-stage optimal policies to be optimal.
\newblock {\em Z. Wahrscheinlichkeitsth}, 32:179--296, 1975.

\bibitem{Sznaier2025DataDrivenLearningControl}
M.~Sznaier, F.~Allgower, A.~C.~B. de~Oliveira, N.~Ozay, and E.~D. Sontag.
\newblock Tutorial: Data driven and learning enabled control.
\newblock In {\em Proceedings of the 64th IEEE Conference on Decision and Control (CDC)}, pages 2858--2873, 2025.

\bibitem{Willard1970GeneralTopology}
Stephen Willard.
\newblock {\em General Topology}.
\newblock Addison-Wesley Publishing Company, Reading, Massachusetts, 1970.

\bibitem{yuksel2023borkar}
S.~Y{\"u}ksel.
\newblock On {B}orkar and {Y}oung relaxed control topologies and continuous dependence of invariant measures on control policy.
\newblock {\em SIAM Journal on Control and Optimization}, 62(4):2367--2386, 2024.

\bibitem{yichen2024}
Y.~Zhou, Y.~Song, and S.~Y\"uksel.
\newblock Robustness to model approximation, empirical model learning, and sample complexity in {W}asserstein regular mdps.
\newblock {\em arXiv preprint arXiv:2410.14116}, 2024.

\end{thebibliography}
\end{document}